\documentclass[11pt]{article}
\usepackage[right=2.5cm,left=2.5cm,top=2cm,bottom=3cm,headsep=2cm,footskip=2cm]{geometry}
\usepackage{graphicx}  
\usepackage{enumerate}
\usepackage{amsmath}
\usepackage{amssymb}
\usepackage{color}
\usepackage{float}
\usepackage{cite}
\usepackage{algorithmic}
\usepackage{algorithm,caption}
\usepackage{mathtools}
\usepackage{hyperref}
\hypersetup{
	colorlinks=true,        % false: boxed links; true: colored links
	linkcolor=black,        % color of internal links
	citecolor=blue,         % color of links to bibliography
	urlcolor=blue           % color of external links
}

\usepackage{authblk}

\title{\bf Zeroth-order optimisation on subsets of symmetric matrices with application to MPC tuning}
\author[1]{Alejandro I. Maass}
\author[1]{Chris Manzie}
\author[2]{Iman Shames}
\author[3]{Hayato Nakada}
\affil[1]{\small Department of Electrical and Electronic Engineering, The University of Melbourne}
\affil[2]{College of Engineering \& Computer Science, The Australian National University}
\affil[3]{Advanced Unit Management System Development Division, Toyota Motor Corporation, Japan}
\affil[*]{e-mails: \small \texttt{\{maassa,manziec\}@unimelb.edu.au;iman.shames@anu.edu.au;hayato\_nakada@mail.toyota.co.jp}  }

\date{}

\newcommand{\qed}{\hfill \ensuremath{\blacksquare}}
\newcommand{\R}{\mathbb{R}}
\newcommand{\N}{\mathbb{N}}

\newcommand{\inner}[1]{\left\langle #1 \right\rangle}

\newcommand{\norm}[1]{\left\|#1\right\|}
\newcommand{\trace}[1]{\mbox{Tr}\left\{#1\right\}}

\newcommand{\diag}[1]{\mbox{diag}\{#1\}}
\newcommand{\Q}{\mathbb{Q}}
\newcommand{\G}{\mathbb{G}}
\renewcommand{\O}{\mathcal{O}}
\renewcommand{\S}{\mathbb{S}}
\renewcommand{\vec}[1]{\mbox{vec}\{#1\}}
\newcommand{\esp}[2]{\mathbb{E}_{#1}\{#2\}}
\newcommand{\espp}[1]{\mathbb{E}\left\{#1\right\}}
\newcommand{\mtxdu}[2]{\begin{bmatrix}#1\\ #2\end{bmatrix}}

\newcommand{\mtxdd}[4]{\begin{bmatrix}#1& #2\\ #3 &#4\end{bmatrix}}

\newcommand{\smallmtx}{\setlength\arraycolsep{2pt}
	\def\arraystretch{1}}

\newtheorem{assum}{Assumption}

\newtheorem{defn}{Definition}
\newtheorem{theo}{Theorem}
\newtheorem{coro}{Corollary}
\newtheorem{lemma}{Lemma}
\newtheorem{remark}{Remark}

\providecommand{\proofname}{Proof}
{%
	\par\noindent{\bfseries\upshape \proofname\ }%
}%
{\qed}

\begin{document}
\maketitle

\begin{abstract}
This paper provides a zeroth-order optimisation framework for non-smooth and possibly non-convex
cost functions with matrix parameters that are real and symmetric.  We provide complexity bounds
on the number of iterations required to ensure a given accuracy level for both the convex and non-convex case. The derived complexity bounds for the convex case are less conservative than available bounds in the literature since we exploit the symmetric structure of the underlying matrix space. Moreover, the non-convex complexity bounds are novel for the class of optimisation problems we consider.   The utility of the framework is evident in the suite of applications that use symmetric matrices as tuning parameters. Of primary interest here  is the challenge of tuning the gain matrices in model predictive controllers, as this is a challenge known to be inhibiting industrial implementation of these architectures. To demonstrate the framework we consider the problem of MIMO diesel air-path control, and consider implementing the framework iteratively ``in-the-loop'' to reduce tracking error on the output channels. Both simulations and experimental results are included to illustrate the effectiveness of the proposed framework over different engine drive cycles.
\end{abstract}

\section{Introduction}
Optimisation with matrix-valued decision variables is a problem that appears in a wide variety of applications. Particularly, this approach can be found in non-negative matrix factorization \cite{lin2007}, low-rank matrix optimisation \cite{zhu2018,zhao2015}, signal processing with independent component analysis \cite{anderson2012}, energy efficient wireless systems \cite{mertikopoulos2017}, graph weight allocation \cite{shafi2011graph,bertrand2013topology}, and controller tuning such as model predictive controllers \cite{butifra12}, linear quadratic regulators \cite{marco2016automatic}, and PIDs for MIMO systems \cite{chaikovskii2009optimal}.
In most of the aforementioned applications, either explicit expressions of the cost function are not available, or derivatives are difficult or infeasible to obtain. 
%For instance, in energy efficient wireless systems, each transmitter controls its individual input signal covariance matrix to minimise some achievable throughput that is often computed from measurements. In graph weight allocation problems and controller tuning problems the underlying dynamics are often too complicated to get an explicit relationship with the design parameters. 
Consequently, we restrict our attention to this class of problems in which the cost function is treated as a black-box. In this case, the cost function directional derivatives can still be approximated by means of cost function values using finite-differences for instance.
Of main interest in this paper is the challenge of tuning the gain matrices in model predictive controllers (MPC), which have the property of being symmetric and positive (semi) definite. Efficient calibration of MPCs is often a difficult task given the large number of tuning parameters and their non-intuitive correlation with the output response.

Existing literature regarding MPC tuning can be broadly divided in heuristic guidelines \cite{garsor10}, analytical studies \cite{shrcoo98,shaebg11,bagkha14}, and meta-heuristic algorithms \cite{satrifa19,yaalod16,jugema14,sukait07,vasyyo08,butifra12,ira2020,sankar2019fast}.
Both heuristic guidelines and analytical studies hold for specific scenarios and MPC formulations, thus trial-and-error approaches are unavoidably adopted in practice. 
Regarding meta-heuristic methods, the works \cite{satrifa19,sukait07,vasyyo08,ira2020,sankar2019fast} focus on set-point tuning based on step response characteristics (e.g. overshoot, settling and rise times, and steady-state errors). 
The authors in \cite{yaalod16} deal with trajectory tuning of MPC in which the trajectories are generated by first-order transfer function models. 
General trajectory tuning is considered in \cite{jugema14} and \cite{butifra12} via particle swarm optimisation (PSO) and gradient-based algorithms, respectively.
The number of particles in PSO algorithms depends on the specific problem and it usually oscillates around 20 to 50 particles. This translates into running 20 to 50 closed-loop experiments per iteration, which is not practical for applications where the plant is in the loop. Gradient-based methods \cite{butifra12} and Bayesian optimisation (BO) methods \cite{lu2020mpc,sorourifar2020data} can also be used as an alternative with less experiments per iteration, but existent results are presented only for vectors of parameters (i.e. diagonal tuning matrices). 
	Lastly, inverse optimal control techniques such as \cite{menner2019constrained} could be used for MPC tuning, see e.g. \cite{ramadan_iMPC}. However, \cite{ramadan_iMPC} does not provide any theoretical guarantees and it relies on a good choice of initial weighting matrices and heuristic guidelines.	
	We note there is currently a lack of a general framework with convergence guarantees that can directly deal with, not only diagonal matrices, but also general subsets of symmetric matrices  in which the algorithm preserves the structure of the tuning matrix at every iteration.

In response to the above discussion, we propose a zeroth-order optimisation framework for non-smooth and possibly non-convex cost functions with matrix parameters that are real and symmetric. We adopt iterates that use zeroth-order matrix oracles to approximate the directional derivative in a randomly chosen direction. This random direction is taken from the Gaussian orthogonal ensemble (GOE), which is the set of matrices with independent normally distributed entries, see e.g. \cite{forrester10,anderson10,tracy96}. 
The algorithm we use is the natural extension of the zeroth-order random search algorithm from \cite{nesspo17} but tailored to deal with cost functions with matrix parameters. For this optimisation framework, we provide convergence guarantees in terms of complexity bounds for both convex and non-convex cost functions.	In \cite{nesspo17}, the counterpart of our optimisation problem is studied, in which functions with vector parameters are used, and the zeroth-order oracle samples from the normal distribution. The authors in \cite{nesspo17} present complexity bounds for both convex and non-convex functions. The bounds for the convex case can be applied to our matrix setup if one vectorises both the parameter and random direction matrices. However, this does not carry the matrix structure and, in fact, these bounds are significantly more conservative than the ones we propose. Moreover, the non-convex bounds proposed in \cite{nesspo17} are only applicable to unconstrained problems and thus cannot be used in our constrained setting. Consequently, we develop novel bounds for the non-convex case, in which the optimisation parameters are constrained to subsets of symmetric matrices. 
	The optimisation framework presented in this paper, when applied to the context of MPC tuning, offers a more general approach with respect to available literature, and it also comes with convergence guarantees. Particularly, it can deal with MPC tuning over trajectories in a general setting as opposed to \cite{satrifa19,sukait07,vasyyo08,ira2020,yaalod16}, and it provides MPC weighting matrices that satisfy the required constraints of being symmetric and positive (semi) definite, and thus provide more degrees of freedom than the usual diagonal choices in \cite{butifra12,jugema14,lu2020mpc,sorourifar2020data}.

The applicability of the proposed approach is illustrated via both simulations and experiments in the problem of diesel air-path control. The MPC parameters are iteratively tuned within the closed-loop setting with the goal of improving the overall tracking performance over an engine drive cycle. 
For the simulation study, we compare the performance of our proposed framework with other gradient-free algorithms available in the literature such as dividing rectangles (DIRECT), PSO, and BO. 
In the experimental testing, we tune MPC controllers in a diesel engine test bench over two segments of the new European driving cycle (NEDC), and one segment of the worldwide harmonised light duty test cycle (WLTC).

A summary of our main contributions can be found below.
\begin{enumerate}
	\item We present theoretical guarantees for a general class of optimisation problems with matrix-valued decision variables on subsets of real and symmetric matrices. The adopted iterates are tailored to this setting in the sense that they produce matrices that belong to the adequate space through projection. When the cost function is convex, the derived complexity bounds are demonstrably less conservative than the ones in \cite{nesspo17} (vector-valued decision variables), since we exploit the symmetric structure of the underlying matrix space. In the non-convex case, we derive new complexity bounds that were not available in the literature for this context.
	\item We illustrate the applicability of our framework by using it to tune the matrix parameters in MPC. This extends current literature such as \cite{satrifa19,sukait07,vasyyo08,ira2020,yaalod16,sankar2019fast}. 
	\item We highlight the efficacy of the approach by a simulation study in which we compare it to other available gradient-free algorithms in the literature (DIRECT, PSO, and BO).
	\item We provide various experimental evaluations in a diesel engine test bench showing significant improvement of MPC tracking performance with only a few iterations of the proposed optimisation approach.
\end{enumerate}

\emph{Organisation:} Section \ref{sec:setup} introduces optimisation framework. The complexity bounds are derived in Section \ref{sec:convergence}. In Section \ref{sec:MPC-diesel}, the proposed optimisation framework is applied to MPC tuning in the context of air-path control. In Section \ref{sec:simulation}, we compare different gradient-free algorithms available in the literature to our framework.
Experimental testing is presented in Section \ref{sec:experiments} for different engine drive cycles. Lastly, conclusions are drawn in Section \ref{sec:conclusion}.

%\subsection{Notation}
\emph{Notation:}
Denote by $\R^{m\times n}$ the set of real matrices of dimension $m\times n$, and $\S^n$ the set of real symmetric matrices of dimension $n\times n$. Let $\N\triangleq \{1,2,\dots\}$ and $\N_0\triangleq \N \cup \{0\}$. 
Given a matrix $M$, $M^\top$ denotes its transpose.
%Given a matrix $M$, $M^\top$ denotes its transpose, and $M\succ 0$ (resp. $M\succeq 0$) denotes that $M$ is positive (resp. semi-) definite. 
We use $\diag{M_1,\dots,M_n}$ to denote the standard block diagonal matrix. The identity and zero matrices of dimension $m\times n$ are denoted by $I_{m\times n}$ and $0_{m\times n}$, respectively.
$\nabla f(X)$ denotes the matrix gradient for a scalar function $f$ with matrix parameter $X$, see e.g. \cite{athsch65,dattorro10}. For matrices $M\in\R^{n\times n}$ and $N\in \R^{n\times n}$, the Frobenius inner product is defined as $\inner{M,N}_F\triangleq \trace{M^\top N}$, which induces the norm $\norm{M}_F\triangleq \sqrt{\inner{M,M}_F}$.
For a random variable $x\in\R$, we write $x\sim \mathcal{N}(0,\sigma^2)$ to say that $x$ is normally distributed with zero mean and variance $\sigma^2$.

\section{Optimisation framework}\label{sec:setup}
We consider the following class of optimisation problems
\begin{align}\label{eq:problem}
	\min_{X\in \Q^n} f(X),
\end{align}
where $X$ is an $n\times n$ matrix parameter, $\Q^n\subset \R^{n\times n}$ is a closed convex set of admissible parameters, and $f$ is a non-smooth and possibly non-convex cost function.
We build upon the random search ideas of \cite{nesspo17}, but tailored to the matrix space $\Q^{n}$. That is, we are interested in solving \eqref{eq:problem} using iterates of the form
\begin{align}\label{eq:iterations}
	X_{k+1} = \pi_{\Q^n}\{ X_k - h_k \O_\mu(X_k,U_k) \},
\end{align}
where $\pi_{\Q^n}$ denotes the Euclidean projection onto the closed convex set $\Q^n$, $h_k$ is a positive scalar known as the step size, and $\O_\mu$ denotes the so-called \emph{zeroth-order random oracle} which is defined as
\begin{align}\label{eq:oracle-matrix}
	\O_\mu(X,U) \triangleq \frac{f(X + \mu U)-f(X)}{\mu}U,
\end{align}
where $\mu>0$ denotes the oracle's precision, and $U$ is a random symmetric matrix that belongs to the Gaussian orthogonal ensemble (GOE) as per Definition \ref{def:GOE} below.
\begin{defn}\label{def:GOE}
	The GOE, denoted by $\G^n$, is the set of random symmetric matrices $U\in\S^n$ with i.i.d. entries such that $[U]_{ii}\sim \mathcal{N}(0,1)$, and $[U]_{ij}\sim \mathcal{N}(0,\frac{1}{2})$ independent of $[U]_{ii}$ for $i < j$, see e.g. \cite{forrester10,anderson10,tracy96}.
\end{defn}

The overall method to solve \eqref{eq:problem} is represented by the pseudo-code above, which we name \emph{zeroth-order random matrix search} (ZO-RMS) algorithm, where $N$ denotes the number of iterations.
Note that the oracle, per iteration, only requires the cost function value at two points instead of first-order or second-order derivative information, and it computes an estimate of the directional derivative in a randomly chosen direction.

\begin{algorithm}[t]
	%\floatname{algorithm}{Zeroth-order Random Matrix Search (ZO-RMS)}
	\caption*{\bf Zeroth-Order Random Matrix Search (ZO-RMS)}
	\label{algorithm1}
	\begin{enumerate}[1:]
		\item Choose $X_0\in \Q^n$, $\mu>0$ and $\{h_k\}_{k\in\N_0}$.
		\item \textbf{for} $k=\{0,\dots,N\}$, $N\in\N$ \textbf{do}
		\item \qquad Generate $U_k\in\G^n$
		\item \qquad $\O_\mu(X_k,U_k) \leftarrow \frac{1}{\mu}(f(X_k+\mu U_k) - f(X_k))U_k$
		\item \qquad $X_{k+1} \leftarrow \pi_{\Q^n}\{X_k - h_k\O_\mu(X_k,U_k)\}$
		\item \textbf{end for}
		\item Return $\hat{X}_N\triangleq \arg\min_{X}[f(X):X\in\{X_0,\dots,X_N\}]$.
	\end{enumerate}
\end{algorithm}

We emphasise that the ZO-RMS algorithm is the natural extension of the random vector search algorithm proposed in \cite{nesspo17} but tailored to matrix decision variables. This extension is essential to our setting since there is currently an absence of a general framework for the constrained class of problems \eqref{eq:problem}, i.e. adequate algorithm and convergence guarantees. In \cite{nesspo17}, the authors provide gradient-free algorithms to solve the optimisation problem $\min_{x\in\bar{\Q}^n\subset\R^n} f(x)$ for both convex and non-convex cost functions, which is the counterpart of \eqref{eq:problem} for vector optimisation variables. A zeroth-order random oracle with a multivariate normally distributed direction $u$ is used. In this paper, we extend the setting in \cite{nesspo17} to optimisation problems over a matrix space using random search matrices $U$ from the GOE, which are the natural counterpart to the normally distributed vector $u$. Note that the complexity bounds provided in \cite{nesspo17} for the convex case are applicable to our setting if we consider $x=\vec{X}$ and $u=\vec{U}$. However, the  matrix structure of $X$ and $U$ is immediately lost. Consequently, the first theoretical goal of this work is to obtain less conservative complexity bounds for the ZO-RMS algorithm when $f$ in \eqref{eq:problem} is convex, by exploiting the structure of the underlying matrix space. With respect to the non-convex case, we note that the complexity bounds available in \cite{nesspo17} are not applicable to our case since they consider an unconstrained non-convex problem in $\R^n$ as opposed to a constrained one such as \eqref{eq:problem}. Therefore, our second theoretical goal is to develop new complexity bounds for the ZO-RMS algorithm when solving constrained non-convex optimisation problems in a matrix space, i.e. problem \eqref{eq:problem} with $f$ non-convex.

	As mentioned in the introduction, the proposed optimisation framework presented here is general and fits many applications. Therefore, the theoretical results presented in the following section will hold for any family of problems that fit the framework. Particularly, in this paper we will illustrate how to apply this framework to tune MPC controllers, since it is an important application for which we can provide interesting contributions validated by simulations and experiments.

\section{Complexity bounds}\label{sec:convergence}
In this section, we study the performance of the ZO-RMS algorithm in terms of complexity bounds that guarantees a given level of accuracy for the algorithm. We provide complexity bounds for both the convex and non-convex case.

Let us first introduce some important definitions.
Definition \ref{def:GOE} implies that the probability distribution $\mathbb{P}(\mathrm{d}U)$ in the GOE $\G^n$ is $\mathbb{P}(\mathrm{d}U) = \frac{1}{\kappa} e^{-\frac{1}{2}\norm{U}_F^2}\mathrm{d}U$,
where $\mathrm{d}U$ is the Lebesgue measure on $\S^n\cong \mathbb{R}^{n(n+1)/2}$, i.e. $\mathrm{d}U \triangleq \prod_{i=1}^{n} \mathrm{d}U_{ii}\prod_{i<j} \mathrm{d}U_{ij}$,
and $\kappa$ is the normalising constant defined as $\kappa \triangleq \int_{\G^n} e^{-\frac{1}{2}\norm{U}_F^2}\mathrm{d}U = 2^{\frac{n}{2}}\pi^{\frac{n^2+n}{4}}$.
Therefore, we can define the expectation of a function $q:\G^n \rightarrow \R$ as
$\mathbb{E}_U\{q(U)\} \triangleq  \int_{\G^n} q(U) \mathbb{P}(\mathrm{d}U) = \frac{1}{\kappa} \int_{\G^n} q(U) e^{-\frac{1}{2}\norm{U}_F^2}\mathrm{d}U$,
and the following moments of interest in $\G^n$,
$m_p \triangleq \mathbb{E}_U\{\norm{U}_F^p\} = \frac{1}{\kappa} \int_{\G^n} \norm{U}_F^p e^{-\frac{1}{2}\norm{U}_F^2}\mathrm{d}U$, for $p\in\N_0$.
We state the following lemma for some moments of interest.
\begin{lemma}\label{lem:hayaka}
	For  every $n\geq 1$, we have that $m_1 \leq \sqrt{\frac{n^2+n}{2}}$, $m_2 = \frac{n^2+n}{2}$, and $m_4=\frac{1}{4}(n^4+2n^3+5n^2+4n)$.
\end{lemma}
\emph{Proof:} 
Define $\psi(p)\triangleq \ln(m_p)$, and note that this function is convex in $p$. Let us write $p=(1-\alpha)\times 0 + \alpha \times 2$ (i.e. $\alpha = p/2$). For $p\in[0,2]$ we have $\alpha\in[0,1]$. Therefore, since $\psi(p)$ is convex, we have that $\psi(p)=\psi(\alpha\times 2 + (1-\alpha)\times 0) \leq \alpha \psi(2) + (1-\alpha)\psi(0)$.
Particularly, $\psi(0) = \ln(m_0)=0$, and $\psi(2)=\ln(m_2)$, then $\psi(p)\leq \alpha \ln(m_2) = \frac{p}{2}\ln(m_2)$, and thus $m_p\leq m_2^{p/2}$. Particularly, $m_1\leq \sqrt{m_2}$. On the other hand, from \cite{hayakawa79} we get that $m_2 = (n^2+n)/2$ and $m_4 = \frac{1}{4}(n^4 + 2n^3 + 5n^2 + 4n)$, which concludes the proof.\qed 

Our analysis relies in the following assumption.
\begin{assum}\label{assu:f}
	\hspace{0cm}
	\begin{enumerate}[(a)]
		\item The set of admissible parameters is a subset of the set of real symmetric matrices of dimension $n\times n$, i.e. $\Q^n\subset \S^n$.
		\item The function $f:\Q^n\rightarrow \R$ is Lipschitz continuous, i.e. $\exists L_0(f)>0$ s.t. $|f(X)-f(Y)|\leq L_0(f)\norm{X-Y}_F$ holds for all $X,Y\in\Q^n$.
	\end{enumerate}
\end{assum}

\subsection{Convex case}\label{sec:convex}

Let $X^*\in \Q^n$ be a stationary point of \eqref{eq:problem}, and $f^*\triangleq f(X^*)$.
In addition, define $U_{0:k}\triangleq \diag{U_0,\dots,U_k}$, a random matrix composed by i.i.d. variables $\{U_k\}_{k\in\N_0}$ associated with each iteration of the scheme.

We are now in a position to state the main result. 
\begin{theo}\label{theo:main-matrix}
	Consider problem \eqref{eq:problem} with $f$ convex.
	Suppose Assumption \ref{assu:f} holds, and let the sequence $\{X_k\}_{k\in\{0,\dots,N\}}$ be generated by the ZO-RMS iterates \eqref{eq:iterations}. Then, for any $N\geq 0$
	\begin{multline}\label{eq:theo6-matrix}
		\tfrac{1}{S_N}\sum_{k=0}^{N} h_k (\phi_k - f^*) \leq \mu L_0(f) \sqrt{\tfrac{n^2+n}{2}}
		+ \tfrac{1}{S_N}\Big( \tfrac{1}{2}\norm{X_0 - X^*}_F^2 \\
		+ \tfrac{1}{8}(n^4+2n^3+5n^2+4n)L_0^2(f)\sum_{k=0}^{N}h_k^2 \Big),
	\end{multline}
	where $S_N\triangleq \sum_{k=0}^{N} h_k$, $\phi_k\triangleq \mathbb{E}_{U_{0:k-1}}\{f(X_k)\}$ for $k\in\N$, and  $\phi_0\triangleq f(X_0)$.
\end{theo}
\emph{Proof:} See Appendix \ref{sec:proof-of-theo}.\qed

A corollary from Theorem \ref{theo:main-matrix} can be obtained which provides expressions for $\mu$, $h_k$, and $N$ that ensure a given level of accuracy for the ZO-RMS algorithm.
\begin{coro}\label{coro:complexity-bound}
	Let $\epsilon>0$ be given.	If $\mu$ and $h_k$ are chosen such that
	\begin{subequations}\label{eq:parameters-matrix}
		\begin{align}
			\mu &\leq \frac{\epsilon}{L_0(f)\sqrt{2(n^2+n)}},\\
			h_k &= \frac{2\bar{r}}{L_0(f)\sqrt{n^4 + 2n^3 + 5n^2 + 4n}\sqrt{N+1}}, 
		\end{align}
	\end{subequations}
	for $k\in\{ 0,\dots,N\}$,
	then, $\mathbb{E}_{U_{0:N-1}}\{ f(\hat{X}_N) \} - f^* \leq \epsilon$ is guaranteed by the ZO-RMS algorithm in
	\begin{align}\label{eq:complexity-bound-me}
		N \geq \frac{L_0^2(f)\bar{r}^2}{\epsilon^2}(n^4 + 2n^3 + 5n^2 + 4n)
	\end{align} 
	iterations, 	where $\bar{r}$ is such that $\norm{X_0 - X^*}_F\leq \bar{r}$.
\end{coro}
\emph{Proof:}
Note that
\begin{align*}
	\mathbb{E}_{U_{0:N-1}}\{ f(\hat{X}_N) \}& - f^* \\
	&\leq \mathbb{E}_{U_{0:N-1}}\Big\{ \tfrac{1}{S_N}\sum_{k=0}^{N} h_k(f(X_k) - f^*) \Big\} \\
	&\stackrel{\eqref{eq:theo6-matrix}}{\leq} \mu L_0(f) \sqrt{\tfrac{n^2+n}{2}} + \tfrac{1}{S_N}\Big( \tfrac{R^2}{2} + \tfrac{1}{8}(n^4 + 2n^3 + 5n^2 + 4n)L_0^2(f)\sum_{k=0}^{N}h_k^2 \Big).
\end{align*}

Suppose the number of steps $N$ is fixed, we can choose $\mu$ and $h_k$ as in \eqref{eq:parameters-matrix}
and further obtain the following bound
\begin{align*}
	\mathbb{E}_{U_{0:N-1}}\{ f(\hat{X}_N) \} - f^* \leq \frac{\epsilon}{2} + \frac{\bar{r}^2}{S_N}  
	= \frac{\epsilon}{2} + \frac{\bar{r}L_0(f)\sqrt{n^4 + 2n^3 + 5n^2 + 4n}}{2\sqrt{N+1}}.
\end{align*}
Therefore, in order to satisfy $\mathbb{E}_{U_{0:N-1}}\{ f(\hat{X}_N) \} - f^* \leq \epsilon$, we need $N \geq  \frac{L_0^2(f)\bar{r}^2}{\epsilon^2}(n^4 + 2n^3 + 5n^2 + 4n) - 1$, concluding the proof.\qed 
%Consequently, picking $N=\tfrac{L_0^2(f)\bar{r}^2}{\epsilon^2}(n^4 + 2n^3 + 5n^2 + 4n)$ satisfies the required accuracy. The proof is now complete.\qed  

The complexity bound \eqref{eq:complexity-bound-me} in Corollary \ref{coro:complexity-bound} corresponds to the number of iterations in which the ZO-RMS algorithm guarantees a given accuracy of $\epsilon$, provided the step size $h_k$ and parameter $\mu$ are chosen as per \eqref{eq:parameters-matrix}.
Certainly, the expressions for $h_k,\mu$, and $N$ depend on the Lipschitz constant $L_0(f)$ and $\bar{r}$ which may be hard to obtain explicitly depending on the application. 
However, these can be numerically bounded  as we explain further below in Section \ref{sec:parameters}.
Since the oracle \eqref{eq:oracle-matrix} is random, the guarantees of the ZO-RMS algorithm hold in expectation. In essence, Corollary \ref{coro:complexity-bound} provides sufficient conditions on the oracle's precision $\mu$ and step size $h_k$ such that we can get sufficiently close to the optimal cost function $f^*$ in $N$ iterations, in average.

\begin{figure}
	\centering 
	\includegraphics[scale=0.8]{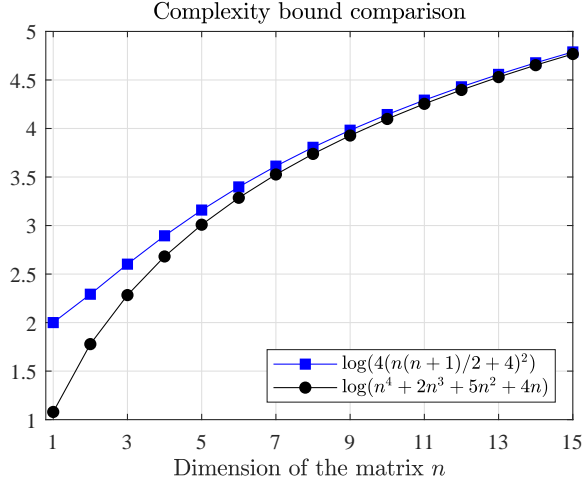}
	\caption{Comparison of the complexity bounds \eqref{eq:complexity-bound-me} and \eqref{eq:nesterov}, see black and blue lines, respectively.}
	\label{fig:comparison}
\end{figure}

We now compare our complexity bound \eqref{eq:complexity-bound-me} to the bounds in \cite{nesspo17} which are applicable to \eqref{eq:problem} (for convex $f$) after grouping the elements of the tuning matrix $X$ into a single vector. Particularly, since $X$ is symmetric, we group the lower triangular elements into a vector of size $n(n+1)/2$ (or often called \emph{half vectorization}), therefore the results in \cite{nesspo17} apply and the following complexity bound holds (cf. Theorem 6 in \cite{nesspo17})
\begin{align}\label{eq:nesterov}
	N_{\text{\tiny \cite{nesspo17}}}\geq \frac{4L_o^2(f)\bar{r}^2}{\epsilon^2}(n(n+1)/2+4)^2.
\end{align}
Fig. \ref{fig:comparison} depicts \eqref{eq:complexity-bound-me} and \eqref{eq:nesterov} for different values of $n$. We can see that our bound \eqref{eq:complexity-bound-me} is less conservative than \eqref{eq:nesterov} for all $n\in\N$, and tends to \eqref{eq:nesterov} as $n\rightarrow\infty$. 
	We emphasise that our optimisation framework is directly tailored to matrix parameters, and the reduction in conservatism with respect to a vector approach such as \cite{nesspo17} follows from exploiting the special structure of the underlying matrix space, namely the symmetry, and the use of iterates that are tailored to it. The reduction in conservativeness is significant, for instance, for $3\times 3$ matrices (i.e. $n=3$), $N_{\text{\tiny \cite{nesspo17}}}/N_{\text{\tiny Cor.\ref{coro:complexity-bound}}}=2.08$. This translates in a sufficient condition that guarantees a level of accuracy in twice less iterations.

\subsection{Non-convex case} 
	Compared to the convex optimisation problem considered in Section \ref{sec:convex}, the analysis for constrained non-convex problems is more challenging. Particularly,  constrained convex problems typically focus on the optimality gap $\espp{f(x)}-f^*$ to measure the convergence rate (as we do in Corollary \ref{coro:complexity-bound}). On the other hand, Nesterov in \cite{nesspo17} provided complexity bounds for  unconstrained non-convex problems using zeroth order iterates, where the optimisation variable is a vector in a subset of $\R^n$. In these unconstrained non-convex problems, the object $\espp{\norm{\nabla f(x)}^2}$ is the typical measure for stationarity. We emphasise that the bounds for unconstrained non-convex problems in \cite{nesspo17} are not applicable to our constrained problem \eqref{eq:problem}, since we need to search over the matrix space $\mathbb{Q}^n$ and use iterates of the form \eqref{eq:iterations} that project onto $\mathbb{Q}^n$. For constrained non-convex problems like ours, a fitting alternative to $\espp{\norm{\nabla f(x)}^2}$ is to consider the so-called \emph{gradient mapping}, see e.g. \cite{liu2018zeroth,ghadimi2016mini}, which is defined as follows
	\begin{align*}
		P_{\mathbb{Q}^n}(X,\mathcal{O}_\mu,h) \triangleq \frac{1}{h}\Big[X - \pi_{\Q^n}\{X - h\mathcal{O}_\mu(X,U)\}\Big],
	\end{align*}
	where $\mathcal{O}_\mu$ is the random zeroth order oracle in \eqref{eq:oracle-matrix}, and recall that $\pi_{\Q^n}$ denotes the Euclidean projection onto the closed convex set $\mathbb{Q}^n$. The natural interpretation of $P_{\mathbb{Q}^n}$ is that it represents the projected gradient, which offers a feasible update from the previous point $X$. The main goal in this section is to provide complexity bounds for the ZO-RMS algorithm \eqref{eq:iterations} in terms of bounding $\espp{\norm{P_{\mathbb{Q}^n}(X,\mathcal{O}_\mu,h)}_F^2}$ when solving the constrained optimisation problem \eqref{eq:problem} with $f$ non-convex.
	
	Before stating our results, we impose an additional assumption.
	\begin{assum}\label{assu:variance}
		$\espp{\norm{\mathcal{O}_\mu(X,U)-\nabla f_\mu(X)}_F^2} \leq \sigma^2$, $\sigma>0$, where $f_\mu$ is the Gaussian approximation of $f$ defined in \eqref{eq:f_mu-matrix}.
	\end{assum}
	Assumption \ref{assu:variance} essentially bounds the variance of the random oracle $\mathcal{O}_\mu$. This assumption is often adopted in non-convex problems, see e.g. \cite{ghadimi2016mini}. If constructing $\sigma$ explicitly is of interest, we can construct it as follows. We note that $\mathcal{O}_\mu$ is an unbiased estimator of $\nabla f_\mu(X)$ (cf. \eqref{eq:(22)-matrix}), i.e. $\espp{\mathcal{O}_\mu(X,U)}=\nabla f_\mu(X)$, then \sloppy  $\espp{\norm{\mathcal{O}_\mu(X,U)-\nabla f_\mu(X)}_F^2} \leq \espp{\norm{\mathcal{O}_\mu(X,U)}_F^2}\stackrel{\text{\tiny Thm.\ref{theo:4-matrix}}}{\leq} (1/4)L_0^2(f)(n^4+2n^3+5n^2+4n)\triangleq \sigma^2$.
	
	The main result for the non-convex case is stated below. It can be seen as the non-convex counterpart of Theorem \ref{theo:main-matrix}.
	\begin{theo}\label{theo:nonconvex}
		Consider problem \eqref{eq:problem} with $f$ non-convex. Suppose Assumptions \ref{assu:f} and \ref{assu:variance} hold, and let the sequence $\{X_k\}_{k\in\{0,\dots,N\}}$ be generated by the ZO-RMS iterates \eqref{eq:iterations}. Then, for any $N\geq 0$,
		\begin{multline}\label{eq:theo-nonconvex}
			\frac{1}{S_N}\sum_{k=0}^{N} h_k \esp{U_{0:k}}{\norm{P_{\mathbb{Q}^n}(X_k,\mathcal{O}_\mu(X_k,U_k),h_k)}_F^2}
			\leq \frac{1}{S_N}\Big[f_\mu(X_0) - f^* + C(\mu)\sum_{k=0}^{N}h_k^2 \Big] + \sigma^2,
		\end{multline}
		where $C(\mu)\triangleq \frac{L_0^3(f)}{4\mu}\sqrt{\frac{n^2+n}{2}}(n^4+2n^3+5n^2+4n)$.
	\end{theo}
	\emph{Proof:} See Appendix \ref{sec:proof-Thm-nonconvex}.\qed 
	
	The next corollary provides a choice of $\mu,h_k$, and $N$ that ensure a given level of accuracy for the ZO-RMS algorithm up to an error of order $\sigma^2$.
	\begin{coro}\label{coro:nonconvex}
		Let $\delta>0$ and $\varepsilon>0$ be given, and $h_k=h$ for all $k\in\{0,\dots,N\}$. If $\mu$ and $h$ are chosen as
		\begin{align}
			\mu &= \frac{\varepsilon}{L_0(f)[(n^2+n)/2]^{1/2}},\label{eq:mu-nonconvex} \\
			h &= \left[ \frac{8 \varepsilon \bar{r}}{(N+1)L_0^3(f)(n^2+n)(n^4+2n^3+5n^2+4n)} \right]^{1/2},\label{eq:h-nonconvex}
		\end{align}
		then $\esp{}{\norm{P_{\mathbb{Q}^n}(X_D,\mathcal{O}_\mu(X_D,U_D),h)}^2} \leq \delta + \sigma^2$ is guaranteed by the ZO-RMS algorithm in 
		\begin{align}\label{eq:N-nonconvex}
			N \geq \frac{L_0^5(f)\bar{r}(n^2+n)(n^4+2n^3+5n^2+4n)}{2\varepsilon \delta^2}
		\end{align}
		iterations, where $$D\triangleq \arg\min_{k\in\{0,\dots,N\}} \norm{P_{\mathbb{Q}^n}(X_k,\mathcal{O}_\mu(X_k,U_k),h)}_F.$$
	\end{coro}
	\emph{Proof:} Recall that $f_\mu$ is a smooth approximation of $f$, and the gap in this approximation can be upper bounded by $|f_\mu(X)-f(X)|\leq \mu L_0(f)\sqrt{\frac{n^2+n}{2}}$, see Lemma \ref{lem:prelminaries-matrix}-(III) in the appendix. Then, to bound this gap by $\varepsilon$ we need to choose $\mu$ as per \eqref{eq:mu-nonconvex}. For this choice of $\mu$, which we denote $\bar{\mu}$, $C(\bar{\mu})=\frac{L_0^4(f)}{4\varepsilon}\left(\frac{n^2+n}{2}\right)(n^4+2n^3+5n^2+4n)$.
	For a constant step size, the right-hand size of \eqref{eq:theo-nonconvex} becomes
	\begin{align*}
		\frac{f_{\bar{\mu}}(X_0)-f^*}{(N+1)h} + C(\bar{\mu})h + \sigma^2 \leq \frac{L_0(f)\bar{r}}{(N+1)h} + C(\bar{\mu})h + \sigma^2.
	\end{align*}
	Let $\rho(h)\triangleq \frac{L_0(f)\bar{r}}{(N+1)h} + C(\bar{\mu})h$. Minimising $\rho(h)$ in $h$ leads to $h^*=\sqrt{L_0(f)\bar{r}/[C(\bar{\mu})(N+1)]}$, which corresponds to \eqref{eq:h-nonconvex}. For this choice of step size, we get $\rho(h^*)=\sqrt{4L_0(f)\bar{r}C(\bar{\mu})/(N+1)}$. Note that $\esp{}{\norm{P_{\mathbb{Q}^n}(X_D,\mathcal{O}_\mu(X_D,U_D),h)}_F^2}\leq \frac{1}{N+1}\sum_{k=0}^{N} \esp{U_{0:k}}{\norm{P_{\mathbb{Q}^n}(X_k,\mathcal{O}_\mu(X_k,U_k),h)}_F^2}\leq \rho(h^*) + \sigma^2$. Then, we can guarantee that $\rho(h^*)\leq \delta$ in $4L_0(f)\bar{r}C(\bar{\mu})/\delta^2$ iterations, which corresponds to \eqref{eq:N-nonconvex}, completing the proof.\qed 
	
	Our analysis shows that ZO-RMS for constrained non-convex problems can suffer an additional error of order $\sigma^2$ which does not appear in the convex case, see Corollary \ref{coro:complexity-bound}. This effect is consistent with the literature on constrained non-convex problems, where this effect has been reported for different algorithms, see e.g. \cite{liu2018zeroth,ghadimi2016mini,liu2018-variance}.

%%%%%%%%%%%%%%%%%%%%%%%%%%%%%%%%%%%%%%%%%%%%%%%%%%%%%%%%%%%%%%%%%%%%%%%%%%%%%%%%%
\section{Application to MPC tuning: air-path control in diesel engines}\label{sec:MPC-diesel}
In this section, we illustrate how to apply the proposed optimisation framework to tune MPC controllers in the context of diesel air-path control.
\subsection{Diesel engine air-path model}
\begin{figure}
	\centering 
	\includegraphics[scale=0.4]{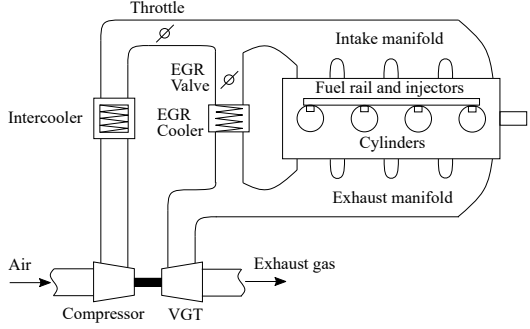}
	\caption{A diesel engine air-path schematic.}
	\label{fig:schematic}
\end{figure}

\begin{figure}
	\centering  
	\includegraphics[scale=0.6]{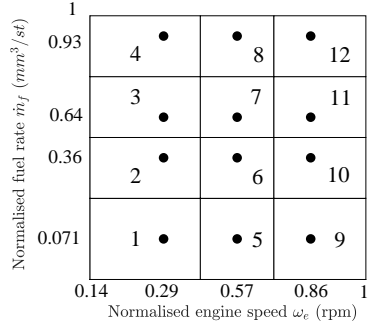}
	\caption{Engine operational space divisions and the corresponding linearisation points.}
	\label{fig:grid}
\end{figure}

A schematic representation of the diesel air-path is shown in Fig. \ref{fig:schematic}. 
The dynamics of a diesel engine air-path are highly non-linear, see e.g. \cite{wahlstrom2011}, and can be captured in the general form
	\begin{subequations}\label{eq:twin}
		\begin{align}
			x_{k+1} &= \mathbf{f}(x_k,u_k,\theta), \\
			y_k &= \mathbf{g}(x_k,u_k,\theta),
		\end{align}	
	\end{subequations}
	where $x_k\in\R^{n_x}, y_k\in \R^{n_y}$, and $u_k\in\R^{n_u}$ are the state, output, and input of the system respectively, at time instant $k\in\N_0$. The engine operational space is typically parametrised by $\theta\triangleq (\omega_e,\dot{m}_f)\in\Theta\subseteq\R^2$, where $\omega_e$ denotes the engine speed, $\dot{m}_f$ denotes the fuel rate, and $\Theta$ is the engine operating space defined as $\Theta\triangleq \{\theta\in\R^2:\omega_{e,\min}\leq \omega_e\leq \omega_{e,\max}, \dot{m}_{f,\min}\leq \dot{m}_f \leq \dot{m}_{f,\max}\}$, for some $\omega_{e,\min},\omega_{e,\max},\dot{m}_{f,\min},\dot{m}_{f,\max}\in\R$. Given these highly non-linear dynamics, a common approach to control-oriented modelling of the air-path is to generate linear models trimmed at various operating points $\theta$. Particularly, we follow \cite{sankar2019fast} and use twelve models to approximate the operating space uniformly. That is, the engine operating space $\Theta$ is divided into twelve regions as per Fig. \ref{fig:grid}, with a linear model representing the engine dynamics in each region. For commercial in confidence purposes, we only show normalised axes in every plot. We emphasise that these regions are chosen to provide adequate coverage over the range of operating points encountered along a drive cycle.
The control-oriented model state consists in the intake manifold pressure $p_{\text{im}}$ (also known as boost pressure), the exhaust manifold pressure $p_{\text{em}}$, the compressor flow $W_{\text{comp}}$, and the EGR rate $y_{\text{EGR}}$. The input consists of the throttle position $u_{\text{thr}}$, EGR valve position $u_{\text{EGR}}$, and VGT valve position $u_{\text{VGT}}$. Lastly, the measured output is $(p_{\text{im}},y_{\text{EGR}})$.

For a given operating point $(\omega_e^\sigma,\dot{m}_f^\sigma)$, $\sigma\in\{1,\dots,12\}$, the engine control unit (ECU) applies certain steady-state actuator values. 
We denote by $\bar{u}^\sigma\in\R^3,\bar{x}^\sigma\in\R^4$, and $\bar{y}^\sigma\in\R^2$ the steady-state values of the input, state, and output respectively at each operating condition $(\omega_e^\sigma,\dot{m}_f^\sigma)$.
Then, by following the system identification procedure detailed in \cite{shsama17}, a linear representation of the non-linear diesel air-path \eqref{eq:twin} trimmed around the grid point $(\omega_e^\sigma,\dot{m}_f^\sigma)$ is given by
\begin{subequations}\label{eq:linear-MVEM}
	\begin{align}
		\tilde{x}_{k+1} &= A^{\sigma}\tilde{x}_k + B^\sigma\tilde{u}_k,\\
		\tilde{y}_k &= C^\sigma\tilde{x}_k + D^\sigma \tilde{u}_k,
	\end{align}
\end{subequations}
where $\sigma\in\{1,\dots,12\}$, $\tilde{x}_k= x_k-\bar{x}^\sigma\in\R^4$ is the perturbed state around the corresponding steady-state $\bar{x}^\sigma$, $\tilde{u}_k= u_k-\bar{u}^\sigma\in\R^3$ is the perturbed input around the corresponding steady-state input $\bar{u}^\sigma$, and $\tilde{y}_k= y_k-\bar{y}^\sigma\in\R^2$ is the perturbed output around the corresponding steady-state output $\bar{y}^\sigma$.


\subsection{MPC formulation}\label{sec:MPC}
For each operating point $(\omega_e^\sigma,\dot{m}_f^\sigma)$, $\sigma\in\{1,\dots,12\}$, and corresponding  model \eqref{eq:linear-MVEM} associated to that region, we formulate an MPC with augmented integrator (see e.g. \cite{rablma09,wang2009model}).
To this end, define the augmented state $\mathbf{x}_k = (\tilde{x}_k,e_k)$, where $e_k$ is the integrator state with dynamics $e_{k+1}= -C^{\sigma}\tilde{x}_k + e_k$.
Define the cost function as
%\begin{multline*}
%J(\tilde{x}_k,\mathbf{u}) \triangleq \tilde{x}_{k+H}^\top P^{\sigma}\tilde{x}_{k+H} \\ + 
%\sum_{i=0}^{H-1} \left( \tilde{x}_{k+i}^\top Q^{\sigma} \tilde{x}_{k+i}+ \tilde{u}_{k+i}^\top R^{\sigma} \tilde{u}_{k+i} \right),
%\end{multline*}
\begin{align*}
	J(\mathbf{x}_k,\mathbf{u}) \triangleq \mathbf{x}_{k+H}^\top {P}^{\sigma}\mathbf{x}_{k+H}  + 
	\sum_{i=0}^{H-1} \left( \mathbf{x}_{k+i}^\top {Q}^{\sigma} \mathbf{x}_{k+i}+ \tilde{u}_{k+i}^\top R^{\sigma} \tilde{u}_{k+i} \right),
\end{align*}
where $H\in\N$ is the prediction horizon, $\mathbf{u}\triangleq \{\tilde{u}_k,\tilde{u}_{k+1},\dots, \tilde{u}_{k+H-1}\}$ is the sequence of control values applied over the horizon $H$, and ${P}^{\sigma}\in\S^{6},{Q}^{\sigma}\in\S^6$, and $R^{\sigma}\in\S^3$ are real symmetric matrices containing the tuning weights. Particularly, we further assume that $R^\sigma$ is positive definite, and that $P^\sigma$ and $Q^\sigma$ are positive semidefinite. 

	The corresponding MPC problem is stated below,
	\begin{align*}
		&\text{minimise}\quad J(\mathbf{x}_k,\mathbf{u})\\
		&\text{subject to}\left\{\begin{array}{l}
			\mathbf{x}_{k+1} = \smallmtx\mtxdd{A^\sigma}{0}{-C^\sigma}{I}\mathbf{x}_k + \mtxdu{B^\sigma}{0}\tilde{u}_k, \\
			\tilde{x}(k) = x(k)-\bar{x}^\sigma,\\
			S_x\mathbf{x}_i + S_u\tilde{u}_i\leq b_i, \forall i\in\{0,\dots,H-1\}, \\
			S_{H}\mathbf{x}_H\leq b_{H},
		\end{array}  \right.
	\end{align*}
	where $S_x\in \R^{18\times 6}$, $S_u\in \R^{18\times 3}$, $b_i\in \R^{18}$, $S_H\in\R^{12\times 6}$, and $b_H\in\R^{12}$ are given by the relevant state and input constraints. The solution to the above optimisation problem yields the optimising control sequence $\mathbf{u}^* =\{ \tilde{u}^*_k,\tilde{u}^*_{k+1}, \dots, \tilde{u}^*_{k+H-1}\}$.
	The first element of the sequence $\mathbf{u^*}$ is applied to the system and the whole process is repeated as $k$ is incremented.

For transient operations between operating points, we use a switching LTI-MPC architecture so that
the controllers are switched based on the current operating condition. As in \cite{sashro19}, the switching LTI-MPC strategy selects the MPC controller at the nearest operating point to the current operating condition.

Under the above scenario, our goal is to use the ZO-RMS algorithm to tune the MPC weighting matrices $\{P^\sigma,Q^\sigma,R^\sigma\}$ in order to achieve a satisfactory tracking performance for a given prediction horizon $H\in\N$. 
Since there are twelve available controllers to tune, we could potentially use the ZO-RMS algorithm to tune all of them. However, since experiments are costly, we will focus on tuning only the controllers that have poor performance given an initial choice of tuning parameters. That is, we utilise an initial calibration denoted by $\{P_0^\sigma,Q_0^\sigma,R_0^\sigma\}$ for every controller $\sigma\in\{1,\dots,12\}$, which we call \emph{baseline controller}. The choice of  $\{P_0^\sigma,Q_0^\sigma,R_0^\sigma\}$ may be based on model dynamics, experience, or using ad-hoc guidelines as per \cite{garsor10}. Then, by looking at a drive cycle response with the baseline controller, we detect the controllers that have poor performance. Let us denote by $\sigma^*$ a controller with unsatisfactory performance and that we attempt to tune with the ZO-RMS algorithm. Therefore, in this context, we want to solve
\begin{align}\label{eq:problem-MPC}
	\min_{P^{\sigma^*},Q^{\sigma^*}\in\S^4_+,R^{\sigma^*}\in\S^3_{++}} g(y(P^{\sigma^*},Q^{\sigma^*},R^{\sigma^*})),
\end{align}
where $g$ is the tracking performance of the MPC defined in \eqref{eq:tracking} below, which depends on a closed-loop response of interest $y$, and $\S^4_+\subset \S^4$ and $\S^3_{++}\subset \S^3$ denote the set of symmetric positive semi-definite matrices and the set of symmetric positive definite matrices, respectively.

As mentioned above, we use the tracking error as the measure of performance for the MPC controller, that is, 
\begin{align}\label{eq:tracking}
	g(y(P^{\sigma^*},Q^{\sigma^*},R^{\sigma^*}))\triangleq  \frac{1}{\sqrt{M}}\sqrt{\sum_{k=0}^{M-1} \left\| y_k(P^{\sigma^*},Q^{\sigma^*},R^{\sigma^*})-y_k^{\text{ref}} \right\|_2^2 }\ ,
\end{align}
where $M\in\N$ is the experiment length, $y_k^{\text{ref}}$ is the vector containing the boost pressure and EGR rate references, respectively, and $y_k(P^{\sigma^*},Q^{\sigma^*},R^{\sigma^*})$ is the process measured output when the $\sigma^*$-th controller is using tuning parameters $P^{\sigma^*},Q^{\sigma^*},R^{\sigma^*}$, and the rest of controllers are using the baseline parameters $\{P_0,Q_0,R_0\}$. We emphasise that for the simulations and experiments below we input normalised signals in \eqref{eq:tracking} so that they are weighted equally.

\subsection{Implementing the ZO-RMS algorithm}
\begin{figure}
	\centering 
	\includegraphics[scale=0.5]{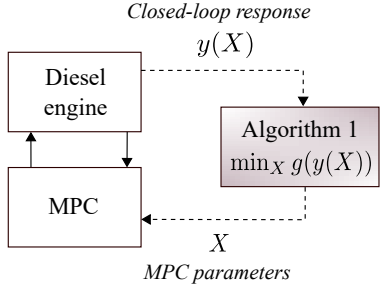}
	\caption{Overall MPC tuning scheme.}
	\label{fig:tuning}
\end{figure}

We now show how to implement the ZO-RMS algorithm for the MPC tuning problem \eqref{eq:problem-MPC}. First note that \eqref{eq:problem-MPC} fits the general optimisation problem \eqref{eq:problem} with $X = \diag{P^{\sigma^*},Q^{\sigma^*},R^{\sigma^*}}$, $f(X)\triangleq g(y(X))$, and 
\begin{align*}
	\Q^n \triangleq \big\{X\in\R^{11\times 11}: X=\diag{P^{\sigma^*},Q^{\sigma^*},R^{\sigma^*}}, P^{\sigma^*},Q^{\sigma^*}\in\S^4_+,R^{\sigma^*}\in\S^3_{++}\big\}.
\end{align*}
Therefore, we can apply ZO-RMS to solve \eqref{eq:problem-MPC}. The overall implementation of ZO-RMS in this context is depicted in Fig. \ref{fig:tuning}. the ZO-RMS algorithm iteratively updates the MPC parameters $\{P^{\sigma^*},Q^{\sigma^*},R^{\sigma^*}\}$ to minimise $f(P^{\sigma^*},Q^{\sigma^*},R^{\sigma^*})$, which is calculated from closed-loop response experiments carried out between the MPC and diesel engine.

It remains to show how to compute the oracle and projection in \eqref{eq:iterations}. Particularly, the oracle in this context is computed as per \eqref{eq:oracle-matrix} with $U = \diag{U^P,U^Q,U^R}$, where $U^P,U^Q\in\G^4$ and $U^R\in\G^3$ (see Definition \ref{def:GOE}). We emphasise that, at each iteration step, the oracle is computed once the entire closed-loop response is available. Particularly, two experiments are required to compute the oracle, one with parameters $\{P_k^{\sigma^*}+\mu U_k^P,Q_k^{\sigma^*}+\mu U_k^Q,R_k^{\sigma^*}+\mu U_k^R\}$ and one with $\{P_k^{\sigma^*},Q_k^{\sigma^*},R_k^{\sigma^*}\}$. Once the two closed-loop responses are finished, then ZO-RMS computes the oracle and next update $\{P_{k+1}^{\sigma^*},Q_{k+1}^{\sigma^*},R_{k+1}^{\sigma^*}\}$.
Two new closed-loop experiments are then performed with these new controller parameters and the process continues iteratively in this fashion.

Since $X$ is block diagonal, the projection onto $\Q^n$ is computed as
$\pi_{\Q^n}\{ X \} = \diag{\pi_{\S^4_+}\{ P \},\pi_{\S^4_+}\{ Q \},\pi_{\S^3_{++}}\{ R\} }$,
where $\pi_{\S^4_+}$ and $\pi_{\S^3_{++}}$ denote the Euclidean projection onto $\S^4_+$ and $\S^3_{++}$, respectively. To compute them, we follow the well known results in \cite{high88}. That is, let $K=V\Lambda V^\top$ be the eigenvalue decomposition of a matrix $K$. Then, $\pi_{\S^+}\{ K \} \triangleq V\max\{0,\Lambda\}V^\top$ and $\pi_{\S^{++}}\{ K \} \triangleq V \max\{d,\Lambda\}V^\top$, $d>0$.

\begin{figure}
	\centering 
	\includegraphics[scale=0.45]{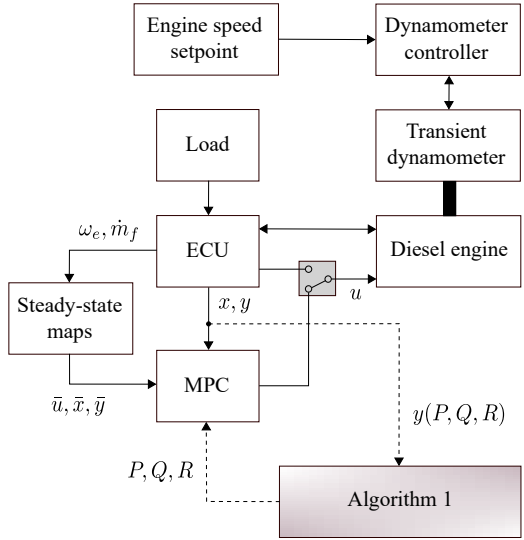}
	\caption{Block diagram of the experimental setting.}
	\label{fig:exp-setup}
\end{figure}

\section{Simulation study}\label{sec:simulation}
We now perform numerical simulations to show the advantage of our proposed MPC tuning framework with respect to other available gradient-free algorithms in the literature. Particularly, we compare our approach to the dividing rectangles (DIRECT) algorithm \cite{jones1993}, particle swarm optimisation algorithms \cite{jugema14}, and Bayesian optimisation algorithms \cite{lu2020mpc,sorourifar2020data}. DIRECT is a sample-based global optimisation method for Lipschitz continuous functions defined over multidimensional domains. It partitions the space into hyperrectangles, each with a sampled point at its centre. The cost function is evaluated at these centrepoints, and then the algorithm keeps sampling and dividing into rectangles until the iteration limit has been reached. PSO algorithms solve the problem by having a population of so-called particles which move along the search space according to the updating rules of their position and velocity. The movement of the particles are guided by their own best known position as well as the entire swarm's best known position in the search space. At each iteration, the cost function is evaluated for every particle in the swarm. 
Lastly, BO seeks to identify the optimal	tuning parameters by strategically exploring and exploiting the parameter space. Exploration aims to evaluate the objective at points in the	decision space with the goal of improving the accuracy of a
	surrogate model of the objective function, while exploitation aims to use the surrogate model to identify decisions that reduce (or increase) the objective function \cite{lu2020mpc}.

Intuition suggests that methods such as DIRECT and PSO would perform many cost function evaluations in order to find the optimal value, since they rely on the number of rectangles/centre points (DIRECT) and particles (PSO). We demonstrate that this is indeed true for the MPC tuning problem described in Section \ref{sec:MPC-diesel}. 

We perform simulations using the tuning scheme from Fig. \ref{fig:tuning}
in which the diesel engine block is simulated via a high-fidelity model, see e.g. \cite{sankar2019fast}. This model is physics-based with the form \eqref{eq:twin}, in which the parameters have been obtained from system identification experiments conducted at Toyota's Higashi-Fuji Technical Centre. The specific equations and parameters are not included due to commercial in confidence purposes.
We focus on tuning controller $\sigma^*=6$ (see Fig. \ref{fig:grid}), over a segment of the NEDC, which is called the urban drive cycle (UDC). We focus on tuning diagonal matrices $\{P^6,Q^6,R^6\}$ with positive   elements so that we can compare to vector-valued methods like DIRECT, PSO, and BO. We emphasise another contribution of the proposed optimisation framework is that it allows the direct tuning of the weighting matrices so that they satisfy the required constraints of symmetricity and positive (semi) definiteness.
For these simulations, we do not include the augmented integrator states. 

The cost function is the tracking error as per \eqref{eq:tracking}, and to assess the complexity of each algorithm we will use the total number of cost function evaluations that takes to achieve the optimal cost function value up to a certain tolerance. The initial calibration is $P_0^6=Q_0^6=\diag{0.01,0,0.2,0.01}$ and $R_0^6=10^{-5}I_{3\times 3}$, and the number of tuning parameters is eleven. 
The ZO-RMS algorithm uses step size $h_k=5\times 10^{-3}/\sqrt{k+1}$ and oracle's precision $\mu = 8\times 10^{-5}$ ($L_0(f)=6,\epsilon=0.008,\bar{r}=2,n=11$).
	For the PSO algorithm we have picked 40 particles as suggested by the literature \cite{jugema14}. For DIRECT, PSO, and BO, we use the decision variable range $[0.01,5]$.

One realisation run of all four algorithms is summarised in Table \ref{tab:algorithm-comparison}, where $f_{\text{opt}}$ denotes the optimal value of the cost function achieved by each algorithm. Next we show the number of cost function evaluations, followed by the total execution time of each algorithm. We can see that in order to achieve a similar value of $f_{\text{opt}}$ (up to a 0.0002 difference), the ZO-RMS algorithm performs 20 function evaluations, compared to 66 for DIRECT, 200 for PSO, and 28 for BO. We emphasise that each function evaluation in this context corresponds to a closed-loop experiment, and thus having to perform many of these is intractable in practice. We can see that BO requires only a few more cost function evaluations than ZO-RMS. It can also be observed that DIRECT, PSO, and BO take longer to execute in comparison to ZO-RMS.

Since PSO, BO, and ZO-RMS are stochastic, we also performed a Monte Carlo simulation of 100 realisations in order to compare the three of them more accurately. These results are listed in Table \ref{tab:montecarlo}. We can see that all achieve a similar optimal value for $f_{\text{opt}}$ in average, but it takes about 199 cost function evaluations in average for PSO to achieve this, and only 22 for BO and 20 for ZO-RMS which makes them more efficient for applications in which the plant is in the loop. 
Overall, BO exhibits comparable performance with respect to ZO-RMS in terms of cost function evaluations and execution time; however, as mentioned in the introduction, these methods only work for vector decision variables and the PSD matrix structure of the MPC tuning matrices is not preserved by the algorithm. They also lack theoretical guarantees in this context.

\begin{table}
	\centering 
	\caption{Comparison between available algorithms.}
	\begin{tabular}{c||c|c|c}
		Algorithm & $f_{\text{opt}}$ & N$^\circ$ of cost fcn. evaluations  & Execution time\\
		\hline 
		DIRECT & 0.0478  & 66 & $\sim$0.6 hours \\
		PSO & 0.0476 & 200 & $\sim$1.6 hours \\
		BO & 0.0477 & 28 & $\sim$0.5 hours\\
		ZO-RMS & 0.0476 & 20 & $\sim$0.3 hours
	\end{tabular}
	\label{tab:algorithm-comparison}
\end{table}

\begin{table}
	\centering 
	\caption{Monte Carlo simulation (100 realisations).}
	\begin{tabular}{c||c|c|c}
		Algorithm & $\esp{}{f_{\text{opt}}}$ & $\textrm{Var}\{f_{\text{opt}}\}$ & $\espp{\text{N$^\circ$ of cost fcn. evaluations}}$ \\
		\hline 
		PSO         & 0.0472 & $5\times 10^{-8}$ & 199  \\
		BO & 0.0477 & $7\times 10^{-8}$ & 22 \\
		ZO-RMS & 0.0496 & $9\times 10^{-6}$ & 20 
	\end{tabular}
	\label{tab:montecarlo}
\end{table}

\section{Experiments}\label{sec:experiments}
The experimental testing of the proposed MPC tuning framework was performed at Toyota's Higashi-Fuji Technical Center in Susono, Japan. The test bench is equipped with a diesel engine, a transient dynamometer, and a dSPACE DS1006 processor board, which is in turn interfaced to the ECU. 
The block diagram in Fig. \ref{fig:exp-setup} depicts the overall experimental setting. 
The ECU logs sensor data from the engine and transmits the current state information to the controller. In addition, the ECU directly controls all engine sub-systems, however, the ECU commands for the three actuators (throttle, EGR valve, and VGT) can be overridden by the MPC commands through a virtual switch in the ControlDesk interface.
The proposed tuning architecture is implemented iteratively in-the-loop as per Fig. \ref{fig:tuning}, in which the diesel engine block is the real engine described above. Particularly, the switched MPC is implemented in real-time on the dSPACE board to generate the required closed-loop drive cycle responses, and the computations in the ZO-RMS algorithm are performed by \textsc{Matlab}, i.e. oracle, projections, and next set of parameters.

\subsection{Choice of algorithm parameters}\label{sec:parameters}
To run ZO-RMS we need to choose the parameters $\mu$ and $h_k$. Corollary \ref{coro:complexity-bound} provides a choice for these parameters so that a given level of accuracy is achieved by the algorithm in $N$ iterations as per \eqref{eq:complexity-bound-me}. However, this choice depends on the Lipschitz constant $L_0(f)$, and the bound $\bar{r}$. These can be numerically bounded depending on the application.
Since we apply the proposed optimisation framework to MPC tuning over a real engine, it is not possible to find an explicit expression for the cost function in terms of $P,Q,R$. Instead, a common approach is to adopt experimental optimisation methods such as the ones in \cite{bunfra16,ahmed2020combining} to ensure that the Lipschitz condition is at least satisfied in the experimental data.
Particularly, since the oracle $\mathcal{O}_\mu$ uses values of the cost function at two different points, we can store these values and use a consistency-check algorithm to see if the Lipschitz condition is verified (see e.g. \cite[Section 4.3]{bunfra16}). This could be done by performing high-fidelity simulations of the closed-loop to compute the oracle at different points, and pick the largest constant that satisfies the Lipschitz condition as a first estimate. Later on in the experiments, we can further adjust this first estimate accordingly, and verify whether the bound is satisfied for the collected experimental measurements. Similarly, $\bar{r}$ can be estimated with the stored data of $X_0$ and $\hat{X}_N$ from multiple simulations so that $\norm{X_0-X^*}\leq \bar{r}$ holds for the experimental space.
In this paper, we initially picked conservative estimates from simulations, and we later refine these heuristically during experimental testing. Online estimation of the Lipschitz constant $L_0(f)$ can improve the aforementioned methods and it is considered as a future direction.

In what follows, we tune three out of the twelve MPC controllers depicted in Fig. \ref{fig:grid} over three different drive cycle segments. Particularly, we tune: controller $\sigma^*=9$ over the middle segment of the UDC (49--117[s]), and we call it UDC2; controller $\sigma^*=5$ over the first 100 seconds of the WLTC, which call WLTC1; and controller $\sigma^*=6$ over the last segment of the UDC (117--195[s]), which we call UDC3.

For all the experiments below we consider an initial calibration $Q_0 = \diag{\bar{Q}_0,\bar{Q}_{e,0}}$, $\bar{Q}_0= \diag{0.01,0.01,0.2,1}$, $\bar{Q}_{e,0}=0.00001\times \diag{1,100}$, $P_0=Q_0,R_0=I_{3\times 3}$ for all twelve controllers, where $\bar{Q}_0$ and $\bar{Q}_{e,0}$ respectively relate to the engine state and integrator state. In addition, all the drive cycle references are generated by the ECU based on a driver model. These models are responsible for converting the drive cycle vehicle speed reference into reference operating points in terms of engine speed and load.

\begin{figure}
	\centering 
	\includegraphics[scale=0.65]{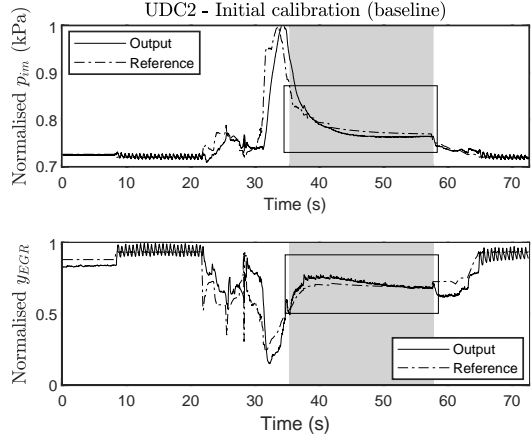}
	\caption{Engine response over the middle UDC segment using the initial parameters $\{P_0,Q_0,R_0\}$ for all twelve controllers.}
	\label{fig:udc2-exp-ini}
\end{figure}

\begin{figure}
	\centering 
	\includegraphics[scale=0.65]{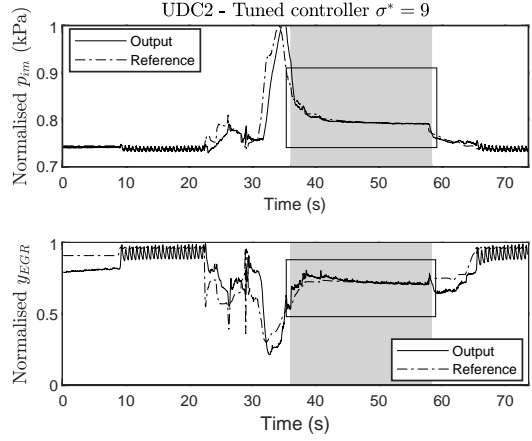}
	\caption{Engine response over the middle UDC segment using the tuned $\{\hat{P}^9,\hat{Q}^9,\hat{R}^9\}$ for controller six, and $\{P_0,Q_0,R_0\}$ for the remaining eleven controllers.}
	\label{fig:udc2-exp-fin}
\end{figure}

%%%%%%%%%%
\subsection{Controller $\sigma^*=9$ over UDC2}
%(2400,10) controller
%$\sigma^*=9$ controller
% N = 11
% f_fin = 3.8207
% f_ini = 4.5595
% f1itr(9) is the minimum and it's 9 since it considers fini. This corresponds to Qopt{8}.
% \mu = 1e-3; hj = 10^(-6)/(sqrt(j+1));
The engine output response to the UDC2 segment with the baseline controller is given in Fig. \ref{fig:udc2-exp-ini}. The grey area illustrates when controller $\sigma^*=9$ is active, and we consider its tracking performance as unsatisfactory (we draw the attention to this with the boxed area). 
%The baseline controller provides a tracking error of $f(P_0,Q_0,R_0)=4.5595$.
We would like to improve the tracking performance for this controller by using ZO-RMS. The parameters used in this experiment are\footnote{We note that for this particular experiment we used a very conservative estimate of the Lipschitz constant based on different experimental tests.} $\mu=0.003\times 10^{-5}$, $h_k=10^{-6}/\sqrt{k+1}$. ($L_0(f)=2.37\times 10^4,\epsilon=0.01,\bar{r}=1.1,n=9$).

In this experiment, we consider $Q^9=\diag{\bar{Q}^9,\bar{Q}_e^9}$ and use ZO-RMS to tune $\{\bar{Q}^9,\bar{Q}_e^9,R^9\}$, whilst $P^9$ is constructed by solving the discrete-time algebraic Riccati equation (DARE) \cite{mayneal00}, $P^9=\texttt{dare}\bigg(\smallmtx\mtxdd{A^9}{0}{-C^9}{I},\mtxdu{B^9}{0},\mtxdd{\bar{Q}^9}{0}{0}{\bar{Q}_e^9},R^9\bigg)$
at each iteration. We emphasise that the proposed tuning framework provides enough flexibility to either tune a single matrix, or all of the MPC matrices depending on the particular application. For instance, for this drive cycle segment, tuning the integrator state matrix $\bar{Q}_e$ helped significantly in improving tracking performance, but for the experiments further below it was not necessary to tune this matrix and it was thus fixed.

The resulting optimal tuning matrices for $N=11$ iterations are given by
\begin{align*}
	\hat{\bar{Q}}^9 &=  \begin{bmatrix}
		0.0611  &  0.0277 &   0.0742 &   0.1389 \\
		0.0277  &  0.1393 &  -0.0554 &   0.0318 \\
		0.0742  & -0.0554 &   0.2177 &   0.0963 \\
		0.1389  &  0.0318 &   0.0963 &   1.0701
	\end{bmatrix},\\
	\hat{R}^9 &= \begin{bmatrix}
		1.8037  & -0.0924  &  0.4081 \\
		-0.0924 &   1.5841 &   0.2252 \\
		0.4081  &  0.2252  &  1.2879
	\end{bmatrix}, \\
	\hat{\bar{Q}}_e^9 &= \begin{bmatrix}
		0.0013  &  0.0014 \\
		0.0014  &  0.0041
	\end{bmatrix}.
	%
	%	\hat{P}^9 &= \texttt{dare}\left(\mtxdd{A^9}{0}{-C^9}{I},\mtxdu{B^9}{0},\mtxdd{\hat{\bar{Q}}^9}{0}{0}{\hat{\bar{Q}}_e^9},\hat{R}^9\right).
\end{align*}

The engine output response over the UDC2 using the above matrices is shown in Fig. \ref{fig:udc2-exp-fin}. The tracking performance has significantly improved in the region where controller nine is acting. Particularly, the tuned matrices provide an improvement of performance\footnote{We compute the percentage of improvement with respect to the baseline controller as: $100\times [f(P_0,Q_0,R_0)-f(\hat{P}^{\sigma^*},\hat{Q}^{\sigma^*},\hat{R}^{\sigma^*})]/f(P_0,Q_0,R_0)$.} of 16.2\% with only eleven iterations.

%%%%%%%%%%%%%%
\subsection{Controller $\sigma^*=5$ over WLTC1}
\begin{figure}
	\centering 
	\includegraphics[scale=0.65]{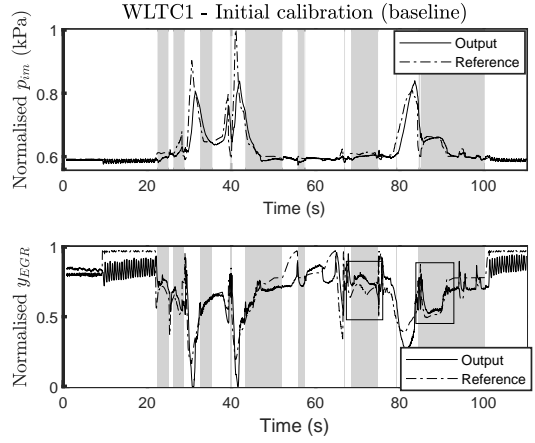}
	\caption{Engine response over the first WLTC segment using the initial parameters $\{P_0,Q_0,R_0\}$ for all twelve controllers.}
	\label{fig:wltc1-exp-ini}
\end{figure}

%(1600,10) controller
%$\sigma^*=5$ controller
% f_ini = 4.1988
% f_fin = 3.2468
% N = 11
% The above function values were computed only in the boxed areas in the picture since one of the non-tuned controllers performed badly in the particular realisation and thus worsen the overall cost function, but actually if this was a simulation that would have not happened and the only thing that actually changes is the grey areas.
% Good picture with better performance is Jitr=8, Qopt{7} 
% 0.34% improvement. But I need to compute it better.
% \mu = 1e-3; hj = 10^(-6)/(sqrt(j+1));
\normalcolor
The engine output response to the WLTC1 segment with the baseline controller is given in Fig. \ref{fig:wltc1-exp-ini}. The grey area illustrates when controller $\sigma^*=5$ is active, and we have boxed the areas with unsatisfactory performance that we would like to improve.
%The baseline controller provides a tracking error of $f(P_0,Q_0,R_0)=4.5595$.
The parameters used in this experiment are $\mu= 2.5\times 10^{-6}$, $h_k=10^{-6}/\sqrt{k+1}$ ($L_0(f)=3450,\epsilon=0.1,\bar{r}=0.1,n=7$).
Similar to the previous experiment, the matrix $P^5$ is constructed by solving the corresponding DARE, and we consider $Q^5=\diag{\bar{Q}^5,\bar{Q}_e^5}$, but now the integrator state matrix $\bar{Q}_e^5$ is not tuned but kept equal to the initial value $\bar{Q}_{e,0}=0.00001\times\diag{1,100}$. 

Consequently, the tuning parameters are $\{\bar{Q}^5,R^5\}$, and their optimal values for $N=11$ iterations are given by 
\begin{align*}
	\hat{\bar{Q}}^5 &= \begin{bmatrix}
		0.0101 &  -0.0078  & -0.0213  &  0.0102 \\
		-0.0078 &   0.0273 &   0.0061 &   0.0170 \\
		-0.0213 &   0.0061 &   0.1760 &   0.0004 \\
		0.0102  &  0.0170  &  0.0004  &  0.9967
	\end{bmatrix}, \\
	\hat{R}^5 &= \begin{bmatrix}
		0.9992  & -0.0066  & -0.0017 \\
		-0.0066 &   0.9839 &   0.0016\\
		-0.0017 &   0.0016 &   0.9806
	\end{bmatrix}.
	%
	%\hat{P}^5 &= \texttt{dare}\left(\mtxdd{A^5}{0}{-C^5}{I},\mtxdu{B^5}{0},\mtxdd{\hat{\bar{Q}}^5}{0}{0}{\bar{Q}_{e,0}},\hat{R}^5\right).
\end{align*}

\begin{figure}
	\centering 
	\includegraphics[scale=0.65]{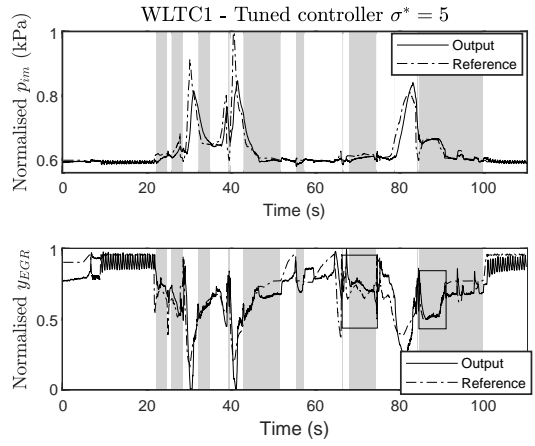}
	\caption{Engine response over the first WLTC segment using the tuned $\{\hat{P}^5,\hat{Q}^5,\hat{R}^5\}$ for controller five, and $\{P_0,Q_0,R_0\}$ for the remaining eleven controllers.}
	\label{fig:wltc1-exp-fin}
\end{figure}

The engine output response over the WLTC1 using the above matrices is shown in Fig. \ref{fig:wltc1-exp-fin}. The tracking performance has improved in the region where controller nine is acting. Particularly, the tuned matrices provide an improvement of performance of 22.67\% with only eleven iterations, as illustrated by the boxed regions in Fig. \ref{fig:wltc1-exp-fin}. We can observe that the tracking has clearly improved in the boxed areas, but we can see that in the grey area 40-50s in Fig. \ref{fig:wltc1-exp-ini} and \ref{fig:wltc1-exp-fin}, the EGR tracking got slightly deteriorated. However, we emphasise that the cost function being minimised is the overall tracking performance considering every grey area where controller $\sigma^*=5$ is acting. This is considered satisfactory since the overall tracking improvement in this case was 22.67\%. Different cost functions or tuning approaches could be potentially used to tackle each region individually.

%%%%%%%%%%%%%%%%%%%%%%%%%%%%%%%%%%%%%%%%%%%%%%%%%%%%%%%%%%%%%%%%%%%%%%%%%%%%%%%%%%%%%%%%%%%
\subsection{Controller $\sigma^*=6$ over UDC3}
%$\sigma^*=6$ controller
% \mu 2.5e-6; hj = 10^{-5}/sqrt(j+1);

\begin{figure}
	\centering 
	\includegraphics[scale=0.65]{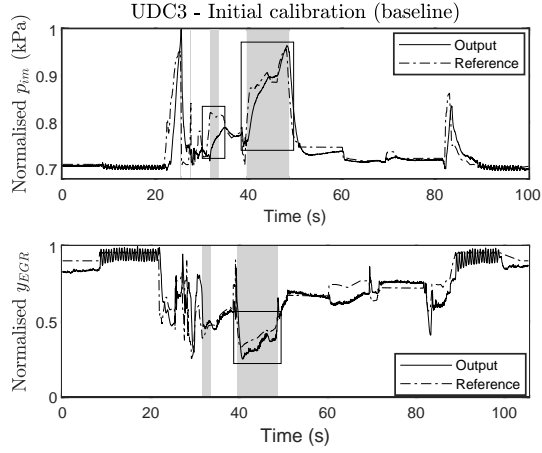}
	\caption{Engine response over the third UDC segment using the initial parameters $\{P_0,Q_0,R_0\}$ for all twelve controllers.}
	\label{fig:udc3-exp-ini}
\end{figure}

The engine output response to the UDC3 segment with the baseline controller is given in Fig. \ref{fig:udc3-exp-ini}. The grey area illustrates when controller $\sigma^*=6$ is active, and we have boxed the areas with unsatisfactory performance. 
%The baseline controller provides a tracking error of $f(P_0,Q_0,R_0)=4.5595$.
The parameters used in this experiment are $\mu=2.5\times 10^{-6}$, $h_k=10^{-5}/\sqrt{k+1}$ ($L_0(f)=2450,\epsilon=0.1,\bar{r}=1.2,n=11$).
The integrator state matrix $\bar{Q}_e$ was not tuned but kept equal to the initial value $\bar{Q}_{e,0}=0.00001\times\diag{1,100}$. As opposed to the previous two experiments, 
we do include $P^5$ in the tuning process for this experiment. Specifically, we construct it as $P^6=\diag{\bar{P}^6,\bar{Q}_{e,0}}$ and tune $\bar{P}^6$.

Consequently, the tuning parameters are \{$\bar{P}^6,\bar{Q}^6,R^6$\}, and their optimal values for $N=8$ iterations are given by 

\begin{align*}
	\hat{\bar{P}}^6 &= \begin{bmatrix}
		0.0502  & -0.0117  & -0.0398  &  0.0057 \\
		-0.0117 &   0.0102 &   0.0146 &   0.0118 \\
		-0.0398 &   0.0146 &   0.2003 &  -0.0426 \\
		0.0057  &  0.0118 &  -0.0426  &  1.0022
	\end{bmatrix},\\
	\hat{\bar{Q}}^6 &= 
	\begin{bmatrix}
		0.4699  &  0.0708 &   0.0192 &   0.0901\\
		0.0708 &   0.2329  & -0.1597 &  -0.1684\\
		0.0192  & -0.1597  &  0.3459  &  0.0016\\
		0.0901 &  -0.1684 &   0.0016 &   1.1911
	\end{bmatrix}, \\
	\hat{R}^6 &= \begin{bmatrix}
		1.1227 &  -0.1419 &   0.0305\\
		-0.1419 &   1.0621 &   0.0534\\
		0.0305 &   0.0534  &  1.2843
	\end{bmatrix}.
\end{align*}

\begin{figure}
	\centering 
	\includegraphics[scale=0.65]{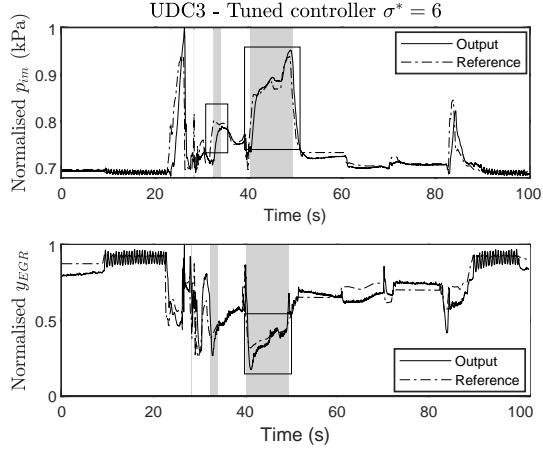}
	\caption{Engine response over the third UDC segment using the tuned $\{\hat{P}^6,\hat{Q}^6,\hat{R}^6\}$ for controller six, and $\{P_0,Q_0,R_0\}$ for the remaining eleven controllers.}
	\label{fig:udc3-exp-fin}
\end{figure}

The engine output response over the UDC3 using the above matrices is shown in Fig. \ref{fig:udc3-exp-fin}. We can see that the tracking performance has improved in the region where controller nine is acting. Particularly, the tuned matrices provide an improvement of performance of 7.73\% with only eight iterations.

The theoretical results in Section \ref{sec:convergence} are sufficient conditions on the step size $h_k$ and oracle's precision $\mu$ that ensure a certain level of accuracy in a fixed number of iterations. We note that the complexity bounds presented in this paper hold for any problem that fits the optimisation framework surrounding \eqref{eq:problem}. Depending on the application, these bounds could become more or less conservative, as we observed in the air-path control case study. In fact, since these are sufficient conditions, we observed that the experimental performance of the algorithm showed tracking error improvement only in a few iterations. Nevertheless, the theoretical results serve as a design tool that guide the choice of algorithm parameters as opposed to a heuristic choice of parameters.

	\begin{remark}
		Overall, we observed that the proposed approach successfully provides improved tracking performance over three different drive cycle segments, i.e. different data sets, which sheds some light on the robustness of the method. However, an interesting future work direction is to provide a theoretical study on robustness guarantees.
	\end{remark}

\section{Conclusion}\label{sec:conclusion}
This paper provided an optimisation framework with theoretical guarantees for the
minimisation of non-smooth and possibly non-convex cost functions with matrix parameters. We applied the proposed algorithm to tune MPCs in the context of air-path control in diesel engines, which is then validated by experimental testing. The algorithm provides improvement of performance with a few iterations, which is demonstrated in different engine drive cycles. Therefore, it creates potential for the development of efficient tuning tools for advanced controllers (and potentially retuning online), even though the theoretical complexity bounds may be large depending on the application.

Future work includes exploring algorithms over other matrix manifolds with practical applications, testing of the proposed algorithms in stochastic drive cycles, robustness guarantees, and online estimation/adaptation of Lipschitz constant $L_0(f)$.

% or


\appendix

\section{Technical results}\label{sec:appendix-preliminaries}
Appendix \ref{sec:appendix-preliminaries} provides the preliminary results needed for Appendices \ref{sec:proof-of-theo} and \ref{sec:proof-Thm-nonconvex}.

Consider a function $f:\Q^n \rightarrow \R$, and let us define the Gaussian approximation of $f$ as
\begin{align}\label{eq:f_mu-matrix}
	f_{\mu}(X) \triangleq  \frac{1}{\kappa} \int_{\G^n} f(X+\mu U)e^{-\frac{1}{2}\norm{U}_F^2}\mathrm{d}U.
\end{align}

We first introduce the following definition.
\begin{defn}\label{def:sub-differential-matrix}
	Let $f$ be a convex function. A matrix $G$ is called a subgradient of function $f$ at $X\in\R^{n\times n}$ if for any $Y\in\R^{n\times n}$ we have	$f(Y) \geq f(X)  + \langle G,Y-X\rangle_F$.
	The set of all subgradients of $f$ at $X\in\R^{n\times n}$ is called subdifferential and is denoted by $\partial f(X)$.
\end{defn}

The following lemma holds for $f_\mu$.
\begin{lemma}\label{lem:prelminaries-matrix}
	\begin{enumerate}[(I)]
		\item If $f$ is convex and $G\in\partial f(X)$, then
		\begin{align}\label{eq:f_mu>f-matrix}
			f_\mu(X) &\geq \frac{1}{\kappa}\int_{\G^n} (f(X)+\mu\langle G,U\rangle_F )e^{-\frac{1}{2}\norm{U}_F^2}\mathrm{d}U = f(X).
		\end{align}
		\item The gradient matrix $\nabla f_\mu$ satisfies
		\begin{align}\label{eq:(22)-matrix}
			\nabla f_\mu ({X})
			= \frac{1}{\kappa} \int_{\G^n} \frac{f({X}+\mu {U}) - f({X})}{\mu}e^{-\frac{1}{2}\norm{{U}}_F^2}{U}\,\mathrm{d}{U}.
		\end{align}
		\item Let $f$ be Lipschitz continuous as per Assumption \ref{assu:f}, then
		\begin{align}\label{eq:(19)-matrix}
			|f_\mu (X) - f(X)| \leq \mu L_0(f) \sqrt{\tfrac{n^2+n}{2}},\quad X\in \Q^n.
		\end{align}
		\item Under Assumption \ref{assu:f}, $f_\mu$ has Lipschitz continuous gradient, i.e. $|f_\mu(Y)-f_\mu(X) - \inner{\nabla f_\mu(X),Y-X}_F| \leq \frac{1}{2} L_1(f_\mu)\norm{X-Y}_F^2$, for all $X,Y\in\Q^n$, with $L_1(f_\mu) = (2L_0(f)/\mu)\sqrt{(n^2+n)/2}$. 
	\end{enumerate}
\end{lemma}
\emph{Proof:} 
\begin{enumerate}[(I)]
	\item We prove this statement as follows.
	\begin{align*}
		f_\mu(X) &\stackrel{\eqref{eq:f_mu-matrix}}{=} \frac{1}{\kappa} \int_{\G^n} f(X+\mu U)e^{-\frac{1}{2}\norm{U}_F^2}\mathrm{d}U \\
		&\hspace{-2.3mm}\stackrel{\mbox{\small Def.\ref{def:sub-differential-matrix}}}{\geq} \frac{1}{\kappa} \int_{\G^n} (f(X)+\langle G,\mu U\rangle_F )e^{-\frac{1}{2}\norm{U}_F^2}\mathrm{d}U \\
		&= f(X) + \frac{1}{\kappa}\int_{\G^n} \langle G,\mu U\rangle_F e^{-\frac{1}{2}\norm{U}_F^2}\mathrm{d}U \\
		&= f(X) + \mu \trace{ \frac{1}{\kappa}\int_{\G^n} G^\top U e^{-\frac{1}{2}\norm{U}_F^2}\mathrm{d}U } \\
		&= f(X),
	\end{align*}
	where the last equality follows from $\esp{U}{U}=0_{n\times n}$.
	\item Define $Y \triangleq X + \mu U$, and note that by the change of variable formula for multivariate integrals, $\mathrm{d}U = (1/\mu^{\bar{n}})\mathrm{d}Y$, where $\bar{n}=n(n+1)/2$. 
	Consequently, 
	\begin{align*}
		f_\mu(X) &= \frac{1}{\kappa} \int_{\G^n} f(Y) e^{-\frac{1}{2\mu^2}\norm{Y-X}_F^2}\times \frac{1}{\mu^{\bar{n}}}\mathrm{d}Y.
	\end{align*}
	Then,
	\begin{align*}
		\nabla f_\mu(X) &= \frac{1}{\mu^{\bar{n}}\kappa}\int_{\G^n} f(Y) e^{-\frac{1}{2\mu^2}\norm{Y-X}_F^2} \times -\frac{1}{2\mu^2}\times \frac{\partial\trace{(Y-X)^2}}{\partial X} \mathrm{d}Y \\
		&= \frac{1}{\mu^{\bar{n}+2}\kappa} \int_{\G^n} f(Y) e^{-\frac{1}{2\mu^2}\norm{Y-X}_F^2} (Y-X) \mathrm{d}Y \\
		&= \frac{1}{\mu \kappa}\int_{\G^n} f(X+\mu U)e^{-\frac{1}{2}\norm{U}_F^2}  U \mathrm{d}U \\
		&= \frac{1}{\kappa} \int_{\G^n} \frac{f(X+\mu U) - f(X)}{\mu}e^{-\frac{1}{2}\norm{U}_F^2}U \mathrm{d} U,
	\end{align*}
	where the last equality follows from $\esp{U}{U}=0_{n\times n}$.
	\item By definition of $f_\mu$ in \eqref{eq:f_mu-matrix} we have
	\begin{align*}
		f_\mu(X) - f(X)  
		= \frac{1}{\kappa}\int_{\G^n} (f(X+\mu U)-f(X))e^{-\frac{1}{2}\norm{U}_F^2}\mathrm{d} U.
	\end{align*}
	Then, by Lipschitz continuity of $f$ (cf. Assumption \ref{assu:f}),
	\begin{align*}
		|f_\mu(X)-f(X)| &\leq \mu L_o(f) \frac{1}{\kappa}\int_{\G^n} \norm{U}_Fe^{-\frac{1}{2}\norm{U}_F^2}\mathrm{d}U \leq \mu L_0(f) \sqrt{\tfrac{n^2+n}{2}},
	\end{align*}
	where the last inequality follows from Lemma \ref{lem:hayaka}.
	\item Indeed, for all $X,Y\in\Q^n$ we have from \eqref{eq:(22)-matrix} and Assumption \ref{assu:f} that
		\begin{align*}
			\norm{\nabla f_\mu(X) - \nabla f_\mu(Y)}_F  \leq \frac{2L_0(f)}{\kappa \mu}\int_{\G^n}\norm{U}_Fe^{-\frac{1}{2}\norm{U}_F^2}\mathrm{d}U ,
		\end{align*}
		and the result follows from Lemma \ref{lem:hayaka}.\qed 
\end{enumerate}

For the random oracle in \eqref{eq:oracle-matrix}, we have the following bound.
\begin{theo}\label{theo:4-matrix}
	If $f$ satisfies Assumption \ref{assu:f}, then 
	\begin{align}\label{eq:(35)-matrix}
		\mathbb{E}_U\{\norm{\O_\mu(X,U)}^2_F\} \leq \frac{1}{4}L_0^2(f)(n^4 + 2n^3 + 5n^2 + 4n).
	\end{align}
\end{theo}
\emph{Proof:} Note that $\norm{\O_\mu(X,U)}_F^2 = \trace{\O_\mu(X,U)^\top \O_\mu(X,U)} = \frac{(f(X+\mu U)-f(X))^2}{\mu^2} \norm{U}_F^2$.
Then, $\mathbb{E}_U\{\norm{\O_\mu(X,U)}^2_F\} \leq \frac{1}{\mu^2} \mathbb{E}_U\{L_0^2(f)\norm{\mu U}_F^2\norm{U}_F^2\} = L_0^2(f) \mathbb{E}_U\{\norm{U}_F^4\}$.
The proof thus follows from Lemma \ref{lem:hayaka}.\qed 

\section{Proof of Theorem 1}\label{sec:proof-of-theo}
We first note that the projection $\pi_{\Q^n}$ in ZO-RMS satisfies
\begin{align}\label{eq:projection-matrix}
	\norm{\pi_{\Q^n}\{X\} - Y}_F\leq \norm{X-Y}_F,
\end{align}
for all $Y\in \Q^n$.

Now, let point $X_k$ with $k\in\N$ be generated after $k$ iterations of ZO-RMS, and define $r_k\triangleq \norm{X_k-X^*}_F$. Then,
\begin{align*}
	r_{k+1}^2 &= \norm{\pi_{\Q^n}\{ X_k-h_k\O_\mu(X_k,U_k) \} - X^*}_F^2 \\
	&\stackrel{\eqref{eq:projection-matrix}}{\leq} \norm{X_k-h_k\O_\mu(X_k,U_k) - X^*}_F^2 \\
	&= r_k^2 - 2h_k \inner{X_k - X^*,\O_\mu(X_k,U_k)}_F  + h_k^2\norm{\O_\mu(X_k,U_k)}_F^2\, .
\end{align*}

Then, taking expectation 
\begin{align*}
	\esp{U_k}{r_{k+1}^2} &\stackrel{\eqref{eq:(35)-matrix}}{\leq} r_k^2 - 2h_k \esp{U_k}{ \inner{X_k - X^*,\O_\mu(X_k,U_k)}_F }  + \tfrac{1}{4}h_k^2L_0^2(f)(n^4 + 2n^3 + 5n^2 + 4n) \\
	&= r_k^2 - 2h_k  \trace{ [X_k-X^*]^\top \esp{U_k}{\O_\mu(X_k,U_k)} }  + \tfrac{1}{4}h_k^2L_0^2(f)(n^4 + 2n^3 + 5n^2 + 4n).
\end{align*}

From \eqref{eq:(22)-matrix}, we get $\esp{U_k}{\O_\mu(X_k,U_k)} = \nabla f_\mu(X_k)$.
Therefore,
\begin{align*}
	\esp{U_k}{r_{k+1}^2} &\leq r_k^2 - 2h_k \inner{X_k-X^*,\nabla f_\mu(X_k)}_F  + \tfrac{1}{4}h_k^2L_0^2(f)(n^4 + 2n^3 + 5n^2 + 4n) \\
	&\stackrel{(a)}{\leq} r_k^2 - 2h_k [f(X_k) - f_{\mu}(X^*)]  + \tfrac{1}{4}h_k^2L_0^2(f)(n^4 + 2n^3 + 5n^2 + 4n),
\end{align*}
where $(a)$ in the last inequality is shown as follows,
\begin{align*}
	f_\mu(X^*) &\stackrel{\text{Def.\ref{def:sub-differential-matrix}}}{\geq} f_\mu(X_k) + \langle \nabla f_\mu(X_k) , X^* - X_k \rangle_F , \\
	&\hspace{0.7mm}\stackrel{\eqref{eq:f_mu>f-matrix}}{\geq}  f(X_k) - \langle X_k - X^* , \nabla f_\mu(X_k) \rangle_F .
\end{align*}

Taking now the expectation in $U_{0:k-1}$, we obtain
\begin{align*}
	\mathbb{E}_{U_{0:k}}\{r_{k+1}^2\} 
	&\leq \mathbb{E}_{U_{0:k-1}}\{ r_k^2 \} - 2h_k (\phi_k-f_\mu(X^*)) + \tfrac{1}{4}h_k^2L_0^2(f)(n^4 + 2n^3 + 5n^2 + 4n).
\end{align*}
Moreover, note that $f_\mu(X^*)\stackrel{\eqref{eq:f_mu>f-matrix}}{\geq} f(X^*)$, thus \eqref{eq:(19)-matrix} implies that $f_\mu(X^*)\leq f^* + \mu L_0(f)\sqrt{(n^2+n)/2}$. Consequently,
\begin{multline}\label{eq:recursion-matrix}
	\mathbb{E}_{U_{0:k}}\{ r_{k+1}^2 \} \leq \mathbb{E}_{U_{0:k-1}}\{ r_k^2 \} - 2h_k (\phi_k-f^*)  
	+ 2h_k\mu L_0(f)\sqrt{\tfrac{n^2+n}{2}}  \\
	+ \tfrac{1}{4}h_k^2L_0^2(f)(n^4 + 2n^3 + 5n^2 + 4n).
\end{multline}

We can now iterate \eqref{eq:recursion-matrix} from $k=0$ to $k=N$ and get
\begin{multline}
	\mathbb{E}_{U_{0:N}}\{ r_{N+1}^2 \} \leq r_0^2 - 2\sum_{k=0}^{N} h_k(\phi_k - f^*)
	  + 2S_N\mu L_0(f)\sqrt{\tfrac{n^2+n}{2}}  \\
	 + \tfrac{1}{4}(n^4 + 2n^3 + 5n^2 + 4n)L_0^2(f) \sum_{k=0}^{N} h_k^2\ . \label{eq:recursion-final-matrix}
\end{multline}
Since $\mathbb{E}_{U_{0:N}}\{ r_{N+1}^2 \}\geq 0$, from \eqref{eq:recursion-final-matrix} we immediately get \eqref{eq:theo6-matrix}, 
%\begin{multline*}
%\tfrac{1}{S_N}\sum_{k=0}^{N} h_k(\phi_k - f^*) 
%%&\leq \mu L_0(f)\sqrt{\tfrac{n^2+n}{2}} - \underbrace{\tfrac{1}{2S_N}\mathbb{E}_{U_{0:N}}(r_{N+1}^2)}_{\geq 0}\\ &\qquad+\tfrac{1}{S_N}\Big(\tfrac{1}{2}\norm{X_0-X^*}_F^2  \\
%%&\hspace{-3mm} + \tfrac{1}{8}(n^4 + 2n^3 + 5n^2 + 4n)L_0^2(f) \sum_{k=0}^{N} h_k^2 \Big)\\
%%%
%\leq \mu L_0(f)\sqrt{\tfrac{n^2+n}{2}} \\
%\qquad +\tfrac{1}{S_N}\Big(\tfrac{1}{2}\norm{X_0-X^*}_F^2 \\
%\hspace{-3mm}  + \tfrac{1}{8}(n^4 + 2n^3 + 5n^2 + 4n)L_0^2(f) \sum_{k=0}^{N} h_k^2 \Big),
%\end{multline*}
completing the proof.\qed

\section{Proof of Theorem \ref{theo:nonconvex}}\label{sec:proof-Thm-nonconvex}
	Define $g_k\triangleq P_{\Q^n}(X_k,\O_\mu(X_k,U_k),h_k)$. From Lemma \ref{lem:prelminaries-matrix}(IV), we have
	\begin{align*}
		f_\mu(X_{k+1}) &\leq f_\mu(X_k) - h_k\inner{\nabla f_\mu(X_k),g_k}_F  + \tfrac{h_k^2 L_1(f_\mu)}{2}\norm{g_k}_F^2.
	\end{align*}
	From non-expansiveness of projection, i.e. \eqref{eq:projection-matrix}, and from the fact that $L_1(f_\mu)=(2L_0(f)/\mu)\sqrt{(n^2+n)/2}$, we have
	\begin{align*}
		f_\mu(X_{k+1}) &\leq f_\mu(X_k) - h_k\inner{\nabla f_\mu(X_k),g_k}_F  + \tfrac{h_k^2 L_0(f)}{\mu}\sqrt{\tfrac{n^2+n}{2}}\norm{\O_\mu(X_k,U_k)}_F^2.
	\end{align*}
	Define $\xi_k\triangleq \O_\mu(X_k,U_k) - \nabla f_\mu(X_k)$, then $\inner{\nabla f_\mu(X_k),g_k}_F=\inner{\O_\mu(X_k,U_k),g_k}_F-\inner{\xi_k,g_k}_F$. Based on \cite[Lemma 1]{ghadimi2016mini} we have $\inner{\O_\mu(X_k,U_k),g_k}\geq \norm{g_k}_F^2 + \frac{1}{h_k}[\textbf{h}(X_{k+1})-\textbf{h}(X_k)]$, where $\textbf{h}(X)=0$ if $X\in\Q^n$ and $\infty$ otherwise.  Since $X_0\in\Q^n$ and based on \eqref{eq:iterations}, we know that $\mathbf{h}(X_k)=0$ for $k\geq 0$. Therefore, 
	\begin{multline}\label{eq:fmu-iter}
		f_\mu(X_{k+1}) \leq f_\mu(X_k)  +h_k\inner{\xi_k,g_k-P_{\Q^n}(X_k,\nabla f_\mu(X_k),h_k)}_F  + h_k\inner{\xi_k,P_{\Q^n}(X_k,\nabla f_\mu(X_k),h_k)}\\ -h_k\norm{g_k}_F^2  + \tfrac{h_k^2 L_0(f)}{\mu}\sqrt{\tfrac{n^2+n}{2}}\norm{\O_\mu(X_k,U_k)}_F^2.
	\end{multline}
	From \cite[Proposition 1]{ghadimi2016mini}, we can write $\inner{\xi_k,g_k-P_{\Q^n}(X_k,\nabla f_\mu(X_k),h_k)}_F\leq \norm{\xi_k}_F^2$. With this fact, we take expectation in \eqref{eq:fmu-iter} and get
	\begin{align}\label{eq:Efmu-iter}
		h_k\esp{U_k}{\norm{g_k}_F^2} \leq 
		f_\mu(X_k)-\esp{U_k}{f_\mu(X_{k+1})} + h_k\sigma^2 + h_k^2 C(\mu),
	\end{align}
	where we also used $\esp{U_k}{\xi_k}=0$, Theorem \ref{theo:4-matrix} to compute $\esp{U_k}{\norm{\O_\mu(X_k,U_k)}_F^2}$, and Assumption \ref{assu:variance} to upper bound $\esp{U_k}{\norm{\xi_k}_F^2}\leq \sigma^2$. Recall that $S_N\triangleq \sum_{k=0}^{N}h_k$ and that $f_\mu(X^*)\geq f^*$. Then, taking expectation in $U_{0:k}$ in \eqref{eq:Efmu-iter} and then taking summation, we obtain \eqref{eq:theo-nonconvex} as required.\qed

%\section{}

\section*{Acknowledgment}
The authors would like to thank the engineering staff at Toyota Higashi-Fuji Technical Centre, Susono, Japan, for their assistance in running the experiments included in this paper.

\bibliography{bibliography}

\begin{thebibliography}{10}

\bibitem{ahmed2020combining}
Mohamed~Osama Ahmed, Sharan Vaswani, and Mark Schmidt.
\newblock Combining bayesian optimization and lipschitz optimization.
\newblock {\em Machine Learning}, 109(1):79--102, 2020.

\bibitem{anderson10}
Greg~W Anderson, Alice Guionnet, and Ofer Zeitouni.
\newblock {\em An introduction to random matrices}, volume 118.
\newblock Cambridge university press, 2010.

\bibitem{anderson2012}
Matthew Anderson, Xi-Lin Li, Pedro Rodriguez, and T{\"u}lay Adali.
\newblock An effective decoupling method for matrix optimization and its
  application to the {ICA} problem.
\newblock In {\em 2012 IEEE International Conference on Acoustics, Speech and
  Signal Processing (ICASSP)}, pages 1885--1888. IEEE, 2012.

\bibitem{athsch65}
M.~Athans and F.C. Schweppe.
\newblock Gradient matrices and matrix calculations.
\newblock Technical report, M.I.T. Lincoln Lab., 1965.

\bibitem{bagkha14}
Peyman Bagheri and Ali Khaki-Sedigh.
\newblock An analytical tuning approach to multivariable model predictive
  controllers.
\newblock {\em Journal of Process Control}, 24(12):41--54, 2014.

\bibitem{bertrand2013topology}
Alexander Bertrand and Marc Moonen.
\newblock Topology-aware distributed adaptation of {L}aplacian weights for
  in-network averaging.
\newblock In {\em 21st European Signal Processing Conference (EUSIPCO 2013)},
  pages 1--5. IEEE, 2013.

\bibitem{bunfra16}
Gene~A Bunin and Gr{\'e}gory Fran{\c{c}}ois.
\newblock Lipschitz constants in experimental optimization.
\newblock {\em arXiv preprint arXiv:1603.07847}, 2016.

\bibitem{butifra12}
Gene~A Bunin, Fernando~Fraire Tirado, Gr{\'e}gory Fran{\c{c}}ois, and Dominique
  Bonvin.
\newblock Run-to-run {MPC} tuning via gradient descent.
\newblock In {\em Computer Aided Chemical Engineering}, volume~30, pages
  927--931. Elsevier, 2012.

\bibitem{chaikovskii2009optimal}
MM~Chaikovskii and Igor'Borisovich Yadykin.
\newblock Optimal tuning of {PID} controllers for {MIMO} bilinear plants.
\newblock {\em Automation and Remote Control}, 70(1):118--132, 2009.

\bibitem{dattorro10}
Jon Dattorro.
\newblock {\em Convex optimization \& Euclidean distance geometry}.
\newblock Meboo Publishing USA, 2011.

\bibitem{forrester10}
Peter~J Forrester.
\newblock {\em Log-gases and random matrices (LMS-34)}.
\newblock Princeton University Press, 2010.

\bibitem{garsor10}
Jorge~L Garriga and Masoud Soroush.
\newblock Model predictive control tuning methods: A review.
\newblock {\em Industrial \& Engineering Chemistry Research}, 49(8):3505--3515,
  2010.

\bibitem{ghadimi2016mini}
Saeed Ghadimi, Guanghui Lan, and Hongchao Zhang.
\newblock Mini-batch stochastic approximation methods for nonconvex stochastic
  composite optimization.
\newblock {\em Mathematical Programming}, 155(1-2):267--305, 2016.

\bibitem{hayakawa79}
T.~Hayakawa and Y.~Kikuchi.
\newblock The moments of a function of traces of a matrix with a multivariate
  symmetric normal distribution.
\newblock {\em South African Statistical Journal}, 13(1):71--82, 1979.

\bibitem{high88}
Nicholas~J Higham.
\newblock Computing a nearest symmetric positive semidefinite matrix.
\newblock {\em Linear algebra and its applications}, 103:103--118, 1988.

\bibitem{ira2020}
Alex~S Ira, Chris Manzie, Iman Shames, Robert Chin, Dragan Ne{\v{s}}i{\'c},
  Hayato Nakada, and Takeshi Sano.
\newblock Tuning of multivariable model predictive controllers through expert
  bandit feedback.
\newblock {\em International Journal of Control}, pages 1--9, 2020.

\bibitem{jones1993}
Donald~R Jones, Cary~D Perttunen, and Bruce~E Stuckman.
\newblock Lipschitzian optimization without the lipschitz constant.
\newblock {\em Journal of optimization Theory and Applications},
  79(1):157--181, 1993.

\bibitem{jugema14}
Gesner A~Nery J{\'u}nior, M{\'a}rcio~AF Martins, and Ricardo Kalid.
\newblock A {PSO}-based optimal tuning strategy for constrained multivariable
  predictive controllers with model uncertainty.
\newblock {\em ISA transactions}, 53(2):560--567, 2014.

\bibitem{lin2007}
Chih-Jen Lin.
\newblock Projected gradient methods for nonnegative matrix factorization.
\newblock {\em Neural computation}, 19(10):2756--2779, 2007.

\bibitem{liu2018-variance}
Sijia Liu, Bhavya Kailkhura, Pin-Yu Chen, Paishun Ting, Shiyu Chang, and Lisa
  Amini.
\newblock Zeroth-order stochastic variance reduction for nonconvex
  optimization.
\newblock {\em Advances in Neural Information Processing Systems},
  31:3727--3737, 2018.

\bibitem{liu2018zeroth}
Sijia Liu, Xingguo Li, Pin-Yu Chen, Jarvis Haupt, and Lisa Amini.
\newblock Zeroth-order stochastic projected gradient descent for nonconvex
  optimization.
\newblock In {\em 2018 IEEE Global Conference on Signal and Information
  Processing (GlobalSIP)}, pages 1179--1183. IEEE, 2018.

\bibitem{lu2020mpc}
Qiugang Lu, Ranjeet Kumar, and Victor~M Zavala.
\newblock {MPC} controller tuning using bayesian optimization techniques.
\newblock {\em arXiv preprint arXiv:2009.14175}, 2020.

\bibitem{marco2016automatic}
Alonso Marco, Philipp Hennig, Jeannette Bohg, Stefan Schaal, and Sebastian
  Trimpe.
\newblock Automatic {LQR} tuning based on gaussian process global optimization.
\newblock In {\em 2016 IEEE international conference on robotics and automation
  (ICRA)}, pages 270--277. IEEE, 2016.

\bibitem{mayneal00}
D.Q. Mayne, J.B. Rawling, C.V. Rao, and P.O.M. Scokaert.
\newblock Constrained model predictive control: Stability and optimality.
\newblock {\em Automatica}, 36:789--814, 2000.

\bibitem{menner2019constrained}
Marcel Menner, Peter Worsnop, and Melanie~N Zeilinger.
\newblock Constrained inverse optimal control with application to a human
  manipulation task.
\newblock {\em IEEE Transactions on Control Systems Technology}, 2019.

\bibitem{mertikopoulos2017}
Panayotis Mertikopoulos, E~Veronica Belmega, Romain Negrel, and Luca
  Sanguinetti.
\newblock Distributed stochastic optimization via matrix exponential learning.
\newblock {\em IEEE Transactions on Signal Processing}, 65(9):2277--2290, 2017.

\bibitem{nesspo17}
Yurii Nesterov and Vladimir Spokoiny.
\newblock Random gradient-free minimization of convex functions.
\newblock {\em Foundations of Computational Mathematics}, 17(2):527--566, 2017.

\bibitem{ramadan_iMPC}
Ahmed Ramadan, Jongeun Choi, Clark~J. Radcliffe, John~M. Popovich, and N.~Peter
  Reeves.
\newblock Inferring control intent during seated balance using inverse model
  predictive control.
\newblock {\em IEEE Robotics and Automation Letters}, 4(2):224--230, 2019.

\bibitem{rablma09}
James~Blake Rawlings and David~Q Mayne.
\newblock {\em Model predictive control: Theory and design}.
\newblock Nob Hill Pub., 2009.

\bibitem{sankar2019fast}
Gokul~S Sankar, Rohan~C Shekhar, Chris Manzie, Takeshi Sano, and Hayato Nakada.
\newblock Fast calibration of a robust model predictive controller for diesel
  engine airpath.
\newblock {\em IEEE Transactions on Control Systems Technology}, 2019.

\bibitem{sashro19}
Gokul~S Sankar, Rohan~C Shekhar, Chris Manzie, Takeshi Sano, and Hayato Nakada.
\newblock Fast calibration of a robust model predictive controller for diesel
  engine airpath.
\newblock {\em IEEE Transactions on Control Systems Technology}, 2019.

\bibitem{satrifa19}
Jose Eduardo~W Santos, Jorge~Ot{\'a}vio Trierweiler, and Marcelo Farenzena.
\newblock Robust tuning for classical {MPC} through the multi-scenarios
  approach.
\newblock {\em Industrial \& Engineering Chemistry Research}, 58(8):3146--3158,
  2019.

\bibitem{shafi2011graph}
S~Yusef Shafi, Murat Arcak, and Laurent El~Ghaoui.
\newblock Graph weight allocation to meet {L}aplacian spectral constraints.
\newblock {\em IEEE Transactions on Automatic Control}, 57(7):1872--1877, 2011.

\bibitem{shaebg11}
Gaurang Shah and Sebastian Engell.
\newblock Tuning {MPC} for desired closed-loop performance for mimo systems.
\newblock In {\em Proceedings of the 2011 American Control Conference}, pages
  4404--4409. IEEE, 2011.

\bibitem{shsama17}
Rohan~C Shekhar, Gokul~S Sankar, Chris Manzie, and Hayato Nakada.
\newblock Efficient calibration of real-time model-based controllers for diesel
  engines--{P}art {I}: Approach and drive cycle results.
\newblock In {\em 2017 IEEE 56th Annual Conference on Decision and Control
  (CDC)}, pages 843--848. IEEE, 2017.

\bibitem{shrcoo98}
Rahul Shridhar and Douglas~J Cooper.
\newblock A tuning strategy for unconstrained multivariable model predictive
  control.
\newblock {\em Industrial \& engineering chemistry research},
  37(10):4003--4016, 1998.

\bibitem{sorourifar2020data}
Farshud Sorourifar, Georgios Makrygirgos, Ali Mesbah, and Joel~A Paulson.
\newblock A data-driven automatic tuning method for mpc under uncertainty using
  constrained bayesian optimization.
\newblock {\em arXiv preprint arXiv:2011.11841}, 2020.

\bibitem{sukait07}
Ryohei Suzuki, Fukiko Kawai, Hideyuki Ito, Chikashi Nakazawa, Yoshikazu
  Fukuyama, and Eitaro Aiyoshi.
\newblock Automatic tuning of model predictive control using particle swarm
  optimization.
\newblock In {\em 2007 IEEE Swarm Intelligence Symposium}, pages 221--226.
  IEEE, 2007.

\bibitem{tracy96}
C.A. Tracy and H.~Widom.
\newblock On orthogonal and symplectic matrix ensembles.
\newblock {\em Communications in Mathematical Physics}, 177(3):727--754, 1996.

\bibitem{vasyyo08}
JH~Van~der Lee, WY~Svrcek, and BR~Young.
\newblock A tuning algorithm for model predictive controllers based on genetic
  algorithms and fuzzy decision making.
\newblock {\em ISA transactions}, 47(1):53--59, 2008.

\bibitem{wahlstrom2011}
Johan Wahlstr{\"o}m and Lars Eriksson.
\newblock Modelling diesel engines with a variable-geometry turbocharger and
  exhaust gas recirculation by optimization of model parameters for capturing
  non-linear system dynamics.
\newblock {\em Proceedings of the Institution of Mechanical Engineers, Part D:
  Journal of Automobile Engineering}, 225(7):960--986, 2011.

\bibitem{wang2009model}
Liuping Wang.
\newblock {\em Model predictive control system design and implementation using
  MATLAB{\textregistered}}.
\newblock Springer Science \& Business Media, 2009.

\bibitem{yaalod16}
Andr{\'e}~Shigueo Yamashita, Paulo~Martin Alexandre, Antonio~Carlos Zanin, and
  Darci Odloak.
\newblock Reference trajectory tuning of model predictive control.
\newblock {\em Control Engineering Practice}, 50:1--11, 2016.

\bibitem{zhao2015}
Tuo Zhao, Zhaoran Wang, and Han Liu.
\newblock A nonconvex optimization framework for low rank matrix estimation.
\newblock In {\em Advances in Neural Information Processing Systems}, pages
  559--567, 2015.

\bibitem{zhu2018}
Zhihui Zhu, Qiuwei Li, Gongguo Tang, and Michael~B Wakin.
\newblock Global optimality in low-rank matrix optimization.
\newblock {\em IEEE Transactions on Signal Processing}, 66(13):3614--3628,
  2018.

\end{thebibliography}
\end{document}